\documentclass[11pt,english]{article}

\usepackage{anysize}
\marginsize{2cm}{2cm}{2.5cm}{2.4cm}

\usepackage[english]{babel}
\usepackage{amsfonts}
\usepackage{latexsym}
\usepackage{amsmath}
\usepackage{mathrsfs}
\usepackage{amssymb}

\let\phi=\varphi
\newcommand{\Z}{{\mathbb Z}}
\newcommand{\R}{{\mathbb R}}
\newcommand{\N}{{\mathbb N}}
\newcommand{\eps}{\varepsilon}
\newcommand{\IP}{{\mathbb P}}

\newcommand{\Po}{{\mathtt P}_{\omega}}

\newcommand{\Poi}{{\mathtt P}_{\omega^{(1)}}}
\newcommand{\Poii}{{\mathtt P}_{\omega^{(2)}}}

\newcommand{\Eo}{{\mathtt E}_{\omega}}
\newcommand{\Eoi}{{\mathtt E}_{\omega^{(1)}}}
\newcommand{\Eoiii}{{\mathtt E}_{\omega^{(3)}}}

\newcommand{\sig}{\sigma}
\newcommand{\gam}{\gamma}

\newcommand{\ga}{\text{\boldmath ${\gamma}$}}

\newcommand{\qed}{\hfill$\Box$\par\medskip\par\relax}

\newcommand{\BB}{{\mathcal B}}

\let\phi=\varphi

\newtheorem{theo}{Theorem}[section]
\newtheorem{lm}{Lemma}[section]

\newtheorem{prop}{Proposition}[section]
\newtheorem{cor}{Corollary}[section]

\allowdisplaybreaks

\title{Random walks with unbounded jumps among random\\ conductances II: Conditional quenched CLT}
\author{Christophe Gallesco \and Serguei Popov}

\begin{document}

\bibliographystyle{plain}

\maketitle
{\footnotesize 

\noindent Department of Statistics, 
Institute of Mathematics, Statistics and Scientific Computation,
University of Campinas--UNICAMP, 
rua S\'ergio Buarque de Holanda 651, 13083--859, Campinas SP,
Brazil\\
\noindent e-mails: \texttt{gallesco@ime.unicamp.br},
\texttt{popov@ime.unicamp.br}

}

\maketitle

\begin{abstract}
We study a one-dimensional random walk among random conductances, 
with unbounded jumps. Assuming the ergodicity of the collection of
conductances and a few other technical conditions (uniform
ellipticity and polynomial bounds on the tails of the jumps)
we prove a quenched \textit{conditional} invariance principle for the
random walk, under the condition that it remains positive until
time~$n$. As a corollary of this result, we
study the effect of conditioning the random walk to exceed level~$n$
before returning to 0 as $n\to \infty$. 
\\[.3cm] \textbf{Keywords:} ergodic environment, unbounded jumps, 
Brownian meander, 3-dimensional Bessel process, hitting
probabilities, crossing time, uniform CLT
\\[.3cm] \textbf{AMS 2000 subject classifications:} 60J10, 60K37
\end{abstract}

\section{Introduction and results}
\label{s_intro}

In this paper, we study one-dimensional random
walks among random conductances, with
unbounded jumps.
This is the continuation of the paper~\cite{GP}, where
we proved a \emph{uniform} quenched invariance principle 
for this model, where ``uniform'' refers to 
the starting position of the walk
(i.e., one obtains the same estimates on the speed of convergence
as long as this position lies in a certain interval around the origin).
 Here, our main results concern the (quenched) limiting law
of the trajectory of the random walk $(X_n, n=0,1,2,\ldots)$ starting
from the origin up to time~$n$, under condition that it remains
positive at the moments $1,\ldots,n$. In Theorem~\ref{Theocond} we
prove that, after suitable rescaling, for a.e.\ environment it
converges to the
\emph{Brownian meander} process, which is, roughly speaking, a
Brownian motion conditioned on staying positive up to some finite
time, and the main result of the paper~\cite{GP} will be an 
important tool for prooving Theorem~\ref{Theocond}. 

This kind of problem was extensively studied for the case of
space-homo\-ge\-neous random walk, i.e., when one can write
$X_n=\xi_1+\cdots+\xi_n$, where the $\xi_i$-s are i.i.d.\ random
variables. These random variables are usually assumed to have
expectation $0$, and to possess some (nice) tail properties. Among
the first papers on the subject we mention~\cite{B72}
and~\cite{Igle}, where the convergence of the rescaled trajectory to
the Brownian meander was proved. Afterwards, finer results (such as
local limit theorems, convergence to other processes if the
original walk is in the domain of attraction of some stable L\'evy
process, etc.) for space-homogeneous random walks were obtained,
see e.g.\ \cite{BD94,C05,CC08,VW09} and references therein. Also,
it is worth noting that in the paper~\cite{B76} the approach
of~\cite{Igle} was substantially simplified by taking advantage of
the homogeneity of the random walk; however, since in our case the
random walk is not space-homogeneous, we rather use methods
similar to those of~\cite{Igle}.

Also, as mentioned in~\cite{GP}, another motivation for this work 
came from Knudsen billiards in
random tubes, see \cite{CP11,CPSV1,CPSV2,CPSV3}.
We refer to Section~1 of~\cite{GP} for the discussion on the relationship
of the present model to random billiards.

Now, we define the model formally.
For $x,y \in \Z$, we denote by $\omega_{x,y}=\omega_{y,x}$ the
conductance
between~$x$ and~$y$. Define $\theta_z\omega_{x,y}=\omega_{x+z,y+z}$, for all $z\in \Z$.
Note that, by Condition~K below, the vectors $\omega_{x,\cdot}$ are
elements of the Polish space $\ell^2(\Z)$.
We assume that $(\omega_{x,\cdot})_{x\in \Z}$ is a stationary 
ergodic (with respect to the family of shifts $\theta$) sequence of
random vectors; $\IP$ stands for the law of this sequence. The collection of all
conductances $\omega=(\omega_{x,y}, x, y\in \Z)$ is called the
\textit{environment}. For all $x\in\Z$, define
$C_x=\sum_{y}\omega_{x,y}$. Given that $C_x<\infty$ for all $x\in\Z$
(which is always so by Condition~K below),
the random walk~$X$ in random environment~$\omega$ is defined through
its transition probabilities
\[
p_{\omega}(x,y)=\frac{\omega_{x,y}}{C_x};
\]
that is, if~$\Po^x$ is the quenched law of the random walk starting
from $x$, we have 
\[
 \Po^x[X_0=x]=1, \quad \Po^x[X_{k+1}=z\mid X_k=y]=p_{\omega}(y,z).
\]
Clearly, this random walk is reversible with the reversible
measure~$(C_x,x\in\Z)$.
 Also, we denote
by~$\Eo^x$ the quenched expectation for the process starting
 from~$x$. When the random walk starts from~$0$, we use shortened
notations $\Po,\Eo$.

In order to prove our results, we need to make two
technical assumptions on the environment:
\medskip

\noindent
{\textbf {Condition~E}.}
There exists $\kappa>0$ such that, $\IP$-a.s., 
$\omega_{0,1}\geq \kappa$. 

\medskip
\noindent{\textbf{Condition~K}.} 
There exist constants $K,\beta>0$ such that $\IP$-a.s.,
$\omega_{0,y}\leq \frac{K}{1+y^{3+\beta}}$,
for all $y\geq 0$.

For future reference, note that combining Conditions E and K we have that there exists $\hat{\kappa}>0$ such that $\IP$-a.s.,
\begin{equation}
\label{Ellip2}
\hat{\kappa} \leq \sum_{y\in \Z} \omega_{0,y}\leq 
\hat{\kappa}^{-1}.
\end{equation}

\medskip

We decided to 
formulate Condition~E this way because, due to the fact that
this work was motivated by random billiards, the main challenge was
to deal with the long-range jumps.
It is plausible that Condition~E could be relaxed to some extent; 
however, for the sake of cleaner presentation of the argument, 
we prefer not trying to deal with \emph{both} long-range jumps and
the lack of nearest-neighbor ellipticity.


Next, for all $n\geq 1$, we define the continuous map $Z^n=(Z^n(t),t\in \R_+)$ as the natural polygonal interpolation of the map $k/n\mapsto
\sig^{-1}n^{-1/2}X_k$ (with $\sig$ from Theorem 1.1 in \cite{GP}). In other words,
\[
\sig \sqrt{n}Z^n_t = X_{\lfloor nt\rfloor}+(nt-\lfloor
nt\rfloor)X_{\lfloor nt\rfloor+1}
\]
with $\lfloor \cdot\rfloor$ the integer part. Also, we denote by~$W$
the standard Brownian motion.

Now, let $\hat{\tau}=\inf\{k\geq 1: X_k\in (-\infty,0]\}$ 
and $ \Lambda_n=\{\hat{\tau}>n\}=\{X_k>0\text{ for all }k=1,\ldots,n\}$.
Consider the conditional quenched probability measure
$Q_{\omega}^n[\;\cdot\;]=\Po[~\cdot\mid \Lambda_n]$, for all
$n\geq 1$.
For each~$n$, the random map $Z^n$ induces a probability measure
$\mu_{\omega}^n$ on $(C[0,1], \BB_1)$, where $\BB_1$
is the Borel $\sig$-algebra on $C[0,1]$ with the supremum norm:
for
any $A\in \BB_1$,
\[
\mu_{\omega}^n(A)=Q_{\omega}^n[Z^n\in A].
\]
Let us next recall the formal definition of the
 Brownian meander $W^+$. For this, let $W$ be a standard Brownian motion and define $\tau_1=\sup\{s\in [0,1]:
W(s)=0\}$ and
$\Delta_1=1-\tau_1$. Then,
\[
W^+(s)=\Delta_1^{-1/2}|W(\tau_1+s\Delta_1)|,\phantom{***}0\leq s\leq
1.
\]
Now, we are ready to formulate the quenched invariance principle for
the random walk conditioned to stay positive, which is the main
result of this paper:
\begin{theo}
\label{Theocond}
Under Conditions~E and~K, we have that, $\IP$-a.s.,
 $\mu_{\omega}^n$ tends weakly to $P_{W^+}$ as $n\to \infty$, where
$P_{W^+}$ is the law of the Brownian meander $W^+$ on $C[0,1]$.
\end{theo}

As a corollary of Theorem~\ref{Theocond}, we obtain a limit
theorem for the process conditioned on crossing a large interval.
Define 
$\hat{\tau}_n=\inf\{k\geq 0: X_k\in [n,\infty)\}$ and $\Lambda'_n=\{\hat{\tau}_n<\hat{\tau}\}$.
We also define $T_n=\inf\{ t>0: Z^{n^2}_t=\sig^{-1}\}$ 
and the stopped process $Y^n_{\cdot}=Z^{n^2}_{\cdot \wedge T_n}$.
Denoting by~$B_3$ the three-dimensional Bessel process (we recall
that $B_3$ is the radial part of a 3-dimensional Brownian motion,
that is, if $(W_1,W_2,W_3)$ is a three-dimensional Brownian motion,
we have $B_3(t)=\sqrt{W_1^2(t)+W_2^2(t)+W_3^2(t)}$) and by 
$\varrho_1=\inf\{t>0: B_3=\sigma^{-1}\}$,
we have
\begin{cor}
\label{Corro}
 Assume Conditions E and K. We have that, $\IP$-a.s., 
under the law $\Po[~\cdot\mid \Lambda_n']$, the couple $(Y^n, T_n)$
converges in law to $(B_3(\cdot \wedge \varrho_1), \varrho_1)$ as
$n\to \infty$. 
\end{cor}

In the next section, we prove some auxiliary results which 
are necessary for the proof of Theorem~\ref{Theocond}. Then, in
Section~\ref{s_proof_Theocond}, we give the
proof of Theorem~\ref{Theocond}. Finally, in
Section~\ref{s_corollary}, we give the proof of
Corollary~\ref{Corro}.

We will denote by $K_1$, $K_2$,~$\dots$ the ``global'' 
constants, that is, those that are used all along the paper and by
$\gam_1$, $\gam_2$,~$\dots$ the ``local" constants, that is, those
that are used only in the subsection in which they appear for the
first time. For the local constants, we restart the numeration in the
beginning of each subsection. 
Besides, to simplify notations, if $x$ is not integer, $\Po^x$
must be understood as $\Po^{\lfloor x\rfloor}$.

\section{Auxiliary results}
\label{s_aux_results}
In this section, we will prove some technical results that 
will be needed later to prove Theorem~\ref{Theocond}.
Let us introduce the following notations. If $A\subset \Z$,
\begin{align}
\label{Temps}
\tau_A=\inf\{n\geq 0: X_n\in A\} \quad \text{and} \quad
\tau^+_A=\inf\{n\geq 1: X_n\in A\}.
\end{align}
Whenever $A=\{x\}$, $x\in \Z$, we write $\tau_x$ 
(respectively, $\tau^+_x$) instead of $\tau_{\{x\}}$ (respectively, $\tau^+_{\{x\}}$).
\medskip

\subsection{Auxiliary environments}
\label{auxenv}
From some fixed environment $\omega$, we are going to introduce three derived environments denoted by $\omega^{(1)}$, $\omega^{(2)}$ and $\omega^{(3)}$ which will be important tools for the proofs of the lemmas in the rest of this section.

Fix two disjoint intervals $B=(-\infty,0]$ and $E=[N, \infty)$ of $\Z$.
For some realization $\omega$ of the environment, consider the new environment~$\omega^{(1)}$ 
obtained from~$\omega$
by deleting all the conductances $\omega_{x,y}$ if~$x$ and~$y$
belong to $(B\setminus\{0\})\cup E$. The reversible measure (up to a constant factor) on this
new environment~$\omega^{(1)}$ is given by
\begin{align*}
C^{(1)}_{0}&=C_{0},\nonumber\\
C^{(1)}_x&=C_x, \phantom{*************}\mbox{if $x\notin B\cup E$,}\nonumber\\
C^{(1)}_x&=\sum_{y\notin (B\setminus\{0\})\cup E}
\omega_{x,y},\phantom{*****}\mbox{otherwise}.
\end{align*}
Now, we define $C^{(1)}_B=\sum_{x\in B}C^{(1)}_x$ and for all $x\in B$,
$\pi_B(x)=C^{(1)}_x/C^{(1)}_B$. Observe that by Conditions~E and~K, $C^{(1)}_B$ is
positive and finite $\IP$-a.s. Hence $\pi_B$ is $\IP$-a.s.\ a
probability measure on $B$.
In the same way we define $\pi_E$ on $E$. For the sake of simplicity
we denote $\Poi^B$ (respectively, $\Poi^E$) instead of
$\Poi^{\pi_B}$ (respectively, $\Poi^{\pi_E}$) for the random walk
on $\omega^{(1)}$ starting with initial distribution~$\pi_B$ (respectively,~$\pi_E$). The same convention will be adopted for environments $\omega^{(2)}$ and $\omega^{(3)}$ defined below.

From the environment $\omega^{(1)}$, we now construct a new environment $\omega^{(2)}$ by setting
if $x > 0$, $y > 0$,
\begin{equation*}
\omega^{(2)}_{x,0}=\sum_{y\in B}\omega^{(1)}_{x,y},
\phantom{*}
\omega^{(2)}_{0,0}=\sum_{y\in
B}\omega^{(1)}_{y,0},
\phantom{*}
\omega^{(2)}_{x,y}=\omega^{(1)}_{x,y}
\end{equation*} 
and $\omega^{(2)}_{x,y}=0$ otherwise. Defining the reversible measure associated to $\omega^{(2)}$ as $C^{(2)}_x=\sum_{y\in \Z}\omega^{(2)}_{x,y}$, for $x\in \Z$, observe in particular that $C^{(2)}_{0}=C^{(1)}_B$ and $C^{(2)}_{x}=C^{(1)}_x$ for $x>0$.

From the environment $\omega^{(1)}$, we finally create a last environment $\omega^{(3)}$ by setting if $x\in
(0,N)$,
\begin{equation*}
\omega^{(3)}_{x,N}=\sum_{y\in E}\omega^{(1)}_{x,y},
\phantom{*}
\omega^{(3)}_{x,0}=\sum_{y\in B}\omega^{(1)}_{x,y}.
\end{equation*} 
Then, let 
\begin{equation*}
\omega^{(3)}_{N,0}=\sum_{y\in
E}\omega^{(1)}_{y,0},
\phantom{*}
\omega^{(3)}_{0,0} =\sum_{y\in
B}\omega^{(1)}_{y,0}.
\end{equation*} 
For $x \in (0,N)$  and $y\in
(0,N)$ we just set $\omega^{(3)}_{x,y}=\omega^{(1)}_{x,y}$ and $\omega^{(3)}_{x,y}=0$ in all other cases. We define the reversible measure associated to $\omega^{(3)}$ as $C^{(3)}_x=\sum_{y\in \Z}\omega^{(3)}_{x,y}$, for $x\in \Z$. Observe in particular that $C^{(3)}_{0}=C^{(1)}_B$, $C^{(3)}_{N}=C^{(1)}_E$ and $C^{(3)}_{x}=C^{(1)}_x$ for $x\in (0, N)$.


\subsection{Crossing probabilities and estimates on the conditional
exit distribution}
\label{s_cross_prob_cond}
Fix $\eps>0$, $n\in \N$ such that $ \eps \sqrt{n}\geq 1$ and take 
$N=\lfloor \eps \sqrt{n} \rfloor$ ($N$ is from section \ref{auxenv}). 
Then define the event
$A_{\eps, n}= \{\tau_{E}<\tau^+_B\}$ ($B$ and $E$ are from section \ref{auxenv}).
For an arbitrary positive integer $M$ define $I_M=[N,N+M]$.
\begin{lm}
\label{SNLcond}
For all $\eta>0$ there exists $M>0$ such that $\IP$-a.s.,
\[
\Po[X_{\tau_{E}}\in I_M\mid A_{\eps,n}]\geq 1-\eta,
\phantom{**}\mbox{for all~$n$ such that $N>1$.}
\]
\end{lm}

\noindent
\textit{Proof.} The proof of this lemma is very similar to the proof of 
Proposition~2.3 of~\cite{GP}. 
Here, we just give the first steps of the proof and then 
indicate the exact place where it matches with the proof of Proposition~2.3 of~\cite{GP}.
First, we write
\begin{align}
\label{XM1cond}
\Po[X_{\tau_{E}}\in I_M\mid A_{\eps,n}]
&=1-\Po[X_{\tau_{E}}\notin I_M\mid A_{\eps,n}]\nonumber\\
&=1-\sum_{y>N+M}\Po[X_{\tau_{E}}=y\mid A_{\eps,n}].
\end{align}
By definition of $\omega^{(1)}$ (cf.\ section \ref{auxenv}), we can couple the random walks in environments $\omega$ and $\omega^{(1)}$  to show that $\Poi[X_{\tau_{E}}=y\mid
A_{\eps,n}]=\Po[X_{\tau_{E}}=y\mid A_{\eps,n}]$. 
Then, by construction of $\omega^{(2)}$, we can couple the random walks in environments $\omega^{(1)}$ and $\omega^{(2)}$  to show that $\Poii[X_{\tau_{E}}=y\mid A_{\eps,n}]=\Poi^{B}[X_{\tau_{E}}=y\mid A_{\eps,n}]$. Thus, we obtain
\begin{align*}
\Poii[X_{\tau_{E}}=y\mid A_{\eps,n}]=\Poi^{B}[X_{\tau_{E}}=y\mid A_{\eps,n}]
&=\sum_{x\in B}\pi_B(x)\Poi^x[X_{\tau_{E}}=y\mid A_{\eps,n}] 
\nonumber\\
&=\sum_{x\in B}\pi_B(x)\Poi^x[X_{\tau_{E}}=y\mid A_{\eps,n}]
\nonumber\\
&\geq \pi_B(0)\Poi[X_{\tau_{E}}=y\mid A_{\eps,n}]\nonumber\\
&= \frac{C_0}{C^{(1)}_B}\Po[X_{\tau_{E}}=y\mid A_{\eps,n}].
\end{align*}
Thus, by~(\ref{XM1cond}) we obtain
\begin{equation*}
\label{PLUC}
 \Po[X_{\tau_{E}}\in I_M\mid A_{\eps,n}]\geq
1-\frac{C^{(1)}_B}{C_0}\sum_{y>N+M}\Poii[X_{\tau_{E}}=y\mid
A_{\eps,n}].
\end{equation*}
Note that, by Condition K and (\ref{Ellip2}), $C^{(1)}_B/C_0\leq \gam_1$ for some constant $\gam_1$. The terms $\Poii[X_{\tau_{E}}=y\mid
A_{\eps,n}]$ can be treated in the 
same way as the terms $\Po^x[X_{\tau_E}=y\mid A_E]$ of equation~(2.6) in \cite{GP}. In particular, following the reasoning anteceding equation (2.9) in \cite{GP}, we can show that 
\begin{equation*}
\Poii[X_{\tau_{E}}=y\mid
A_{\eps,n}]=\frac{C^{(1)}_y\Poii^y[\tau_{0}<\tau_E^+]}{C^{(1)}_E\Poii^E[\tau_{0}<\tau_E^+]}.
\end{equation*}
Then, the numerator and denominator of the above equation can be treated by using the same techniques as those used to treat~(2.9) in \cite{GP}. 
%
\qed
\begin{lm}
\label{Aen}
There exists a positive constant~$K_1$ such that we have, $\IP$-a.s., 
$\Po[A_{\eps, n}] \geq K_1N^{-1}$
for all~$n$ such that $N>1$.
\end{lm}

\noindent
\textit{Proof.} Recall that $\Po[A_{\eps,n}]=\Po[\tau_E<\tau_B^+]$. We can couple the random walks in environments~$\omega$ 
and~$\omega^{(1)}$ (cf.\ section \ref{auxenv}) to show that $\Po[\tau_E<\tau_B^+]=\Poi[\tau_E<\tau_B^+]$.

Let us denote by $\Gamma_{z',z''}$ the set of finite paths 
$(z',z_1,\ldots,z_k,z'')$ such that $z_i\notin B\cup E\cup\{z',z''\}$
for all $i=1,\ldots,k$.
Let $\ga=(z',z_1,\dots,z_k,z'')\in \Gamma_{z',z''}$ and
define 
\[
 \Poi^{z'}[\ga]:=\Poi^{z'}[X_1=z_1,\dots,X_k=z_k,X_{k+1}=z''].
\]
By reversibility we obtain
\begin{align}
\label{XOR}
\Poi[\tau_E<\tau_B^+]
&=\sum_{z\in E}\sum_{\ga\in \Gamma_{0,z}}\Poi[\ga] \nonumber\\
&=\sum_{z\in E} \sum_{\ga\in \Gamma_{z,0}}\frac{C^{(1)}_z}{C_0}
\Poi^z[\ga] \nonumber\\
&=\frac{C^{(1)}_E}{C_0}\sum_{z\in E} \pi_E(z)
\sum_{\ga\in \Gamma_{z,0}}\Poi^z[\ga] \nonumber\\
&=\frac{C^{(1)}_E}{C_0}\Poi^{E}[\tau_B<\tau_E^+, X_{\tau_B}=0].
\end{align}
Now, define $B'=(-\infty,1]$. 
We have
\begin{align}
\label{WART}
\Poi^{E}[\tau_B<\tau^+_E,X_{\tau_B}=0]
&=\Poi^{E}[\tau_B<\tau^+_E,
\tau_{B'}<\tau^+_E,X_{\tau_B}=0]\nonumber\\
&=\Poi^{E}[\tau_{B'}<\tau^+_E] 
\Poi^{E}[\tau_B<\tau^+_E,X_{\tau_B}=0\mid \tau_{B'}<\tau^+_E]\nonumber\\
&\geq \Poi^{E}[\tau_{B}<\tau^+_E] 
\Poi^{E}[\tau_B<\tau^+_E,X_{\tau_B}=0\mid \tau_{B'}<\tau^+_E].
\end{align}
Let us treat the term $\Poi^{E}[\tau_{B}<\tau^+_E]$. By definition of $\omega^{(3)}$ (cf.\ section \ref{auxenv}), we can couple the random walks in environments $\omega^{(1)}$ and $\omega^{(3)}$  to show that $\Poi^{E}[\tau_{B}<\tau^+_E]=\mathtt{P}_{{\omega}^{{(3)}}}^{N}[\tau_{0}<\tau^+_{N}]$. We obtain
\begin{equation}
C^{(1)}_{E}\Poi^{E}[\tau_{B}<\tau^+_E]=C^{(1)}_{E}\mathtt{P}_{{\omega}^{{(3)}}}^{N}[\tau_{0}<\tau^+_{N}]=C^{(3)}_{N}\mathtt{P}_{{\omega}^{{(3)}}}^{N}[\tau_{0}<\tau^+_{N}]=C_{\text{eff}}(1,N) 
\end{equation}
where $C_{\text{eff}}(1,N)$ is the effective conductance between the points $1$ and $N$ of the electrical network associated to ${\omega}^{(3)}$ (cf. \cite{Dosnel}, section 3.4).
Using Condition~E, we obtain
\begin{equation*}
C_{\text{eff}}(1,N) \geq 
\Big(\sum_{i=1}^{N-1}\omega_{i,i+1}^{-1}\Big)^{-1} \geq
\frac{\kappa}{N-1}.
\end{equation*}
Therefore, there exists a constant~$\gam_1$ such that, whenever $N>1$
\begin{equation}
\label{condeffec}
C^{(1)}_{E}\Poi^{E}[\tau_{B'}<\tau^+_E] \geq \frac{\gam_1}{N}.
\end{equation}
Let us treat the term 
$\Poi^{E}[\tau_B<\tau^+_E,X_{\tau_B}=0\mid \tau_{B'}<\tau^+_E]$.
We have by the Markov property 
\begin{align*}
\lefteqn{\Poi^{E}[\tau_B<\tau^+_E,X_{\tau_B}=0 \mid
\tau_{B'}<\tau^+_E]}\phantom{********}\nonumber\\
&=\sum_{y\in \{0,1\}}\Poi^{E}[\tau_B<\tau^+_E, X_{\tau_B}=0,
X_{\tau_{B'}}=y\mid \tau_{B'}<\tau^+_E]\nonumber\\
&=\sum_{y\in\{0,1\}}\Poi^{E}[\tau_B<\tau^+_E, X_{\tau_B}=0,\mid 
X_{\tau_{B'}}=y, \tau_{B'}<\tau^+_E]\Poi^{E}[X_{\tau_{B'}}=y\mid
\tau_{B'}<\tau^+_E]\nonumber\\
&=\sum_{y\in\{0,1\}}\Po^{y}[\tau_B<\tau_E,
X_{\tau_B}=0]\Poi^{E}[X_{\tau_{B'}}=y\mid
\tau_{B'}<\tau^+_E]\nonumber\\
&\geq \min_{y\in\{0,1\}}
\Po^{y}[\tau_B<\tau_E,X_{\tau_B}=0]\nonumber\\
&\geq \Po^1[X_1=0].
\end{align*}
By Condition~E and (\ref{Ellip2}), this last probability is bounded from below by 
the constant~$\kappa \hat{\kappa}$.
Thus, combining this last result with~(\ref{XOR}), (\ref{WART}), (\ref{condeffec})
and, since by (\ref{Ellip2}) we have $C_0\leq \hat{\kappa}^{-1}$,
it follows that $\IP$-a.s.,
\begin{equation*}
\Po[A_{\eps, n}] \geq \frac{\gam_1\kappa \hat{\kappa}^2 }{N}.
\end{equation*}
This concludes the proof of Lemma~\ref{Aen}. 
\qed

\begin{lm}
\label{Little}
There exists a positive constant~$K_2$ such that we have, $\IP$-a.s., 
\[
\Eo[\tau_B^+\wedge \tau_E]\leq K_2N
\]
for all $n$ such that $N>1$.
\end{lm}

\noindent
\textit{Proof.} First notice that by construction of $\omega^{(1)}$ (cf.\ section \ref{auxenv}), we can couple the random walks in environments $\omega$ and $\omega^{(1)}$ to show that $\Eo[\tau_B^+\wedge \tau_E]=\Eoi[\tau_B^+\wedge \tau_E]$. Hence, we obtain
\begin{align*}
\Eoi^{B}[\tau_B^+\wedge \tau_E]
&= \sum_{y\in B} \pi_B(y)\Eoi^y[\tau_B^+\wedge \tau_E]\nonumber\\
&=\pi_B(0)\Eoi^0[\tau_B^+\wedge \tau_E] + \sum_{y\in B\setminus
\{0\}} \pi_B(y)\Eoi^y[\tau_B^+\wedge \tau_E]\nonumber\\
&\geq \pi_B(0)\Eo[\tau_B^+\wedge \tau_E].
\end{align*}
Therefore, we obtain
\begin{equation}
\label{LI1}
\Eo[\tau_B^+\wedge \tau_E] \leq \frac{\Eoi^{B}[\tau_B^+\wedge
\tau_E]}{\pi_B(0)}.
\end{equation}
Then, observe that
\begin{equation}
\label{LI2}
\Eoi^{B}[\tau_B^+\wedge \tau_E]
= \Eoiii[\tau_{0}^+\wedge\tau_{N}].
\end{equation}


 We are going to bound the right-hand side term of (\ref{LI2}) from above. Before this, we make a brief digression to study the invariant measure of a particular process of interest. 
 
 Consider the following particle system in continuous time on the interval $[0,N]$ of $\Z$. Suppose that
we have injection (according to some Poisson process) and absorption
of particles at states~$0$ and~$N$. Once injected,
particles move according to transition rates given by $q_{x,y}=\omega^{(3)}_{x,y}/C^{(3)}_x$, for $(x,y)\in \{0,\dots,N\}^2$, until they reach~$0$ or~$N$. We suppose that injections at $0$ and~$N$ happen
accordingly to independent Poisson processes with rates
respectively $\lambda_{0}=C^{(3)}_0$ and $\lambda_{N}=C^{(3)}_N$.
We are interested in the continuous time Markov process 
$\big(\eta(t)=((\eta_{0}(t),\dots, \eta_{N}(t)), 
t\geq 0\big)$
with state space $\Omega=\Z_+^{\{0,\dots,N\}}$ where $\eta_i(t)$ represents the number of particles in $i$ at time $t$. 
Hereafter, for $(i,j) \in \{0,\dots, N\}^2$, we will use the symbol $\eta^{i,j}$ to denote the configuration obtained from $\eta$ by moving a particle from site $i$ to site $j$, i.e., if for example $i<j$,
$\eta^{i,j} = (\eta_0,\dots , \eta_i-1,\dots, \eta_{j}+1,\dots,\eta_{N})$. We also define $\eta^{i,+} = (\eta_0,\dots, \eta_i+1, \dots, \eta_{N})$ and $\eta^{i,-} = (\eta_0, \dots,\eta_i-1,\dots, \eta_{N})$ for $i\in \{0,\dots, N\}$.
The generator of this process defined by its action on functions $f: \Omega\to \R$ is given by
\begin{align}
\mathcal{L}f(\eta)
&=\lambda_{0}(f(\eta^{0,+})-f(\eta))+\sum_{i=0}^{N}\eta_iq_{i,0}(f(\eta^{i,-})-f(\eta))\nonumber\\
&\phantom{**}+\sum_{i=1}^{N-1}\sum_{j=1}^{N-1}\eta_jq_{j,i}(f(\eta^{j,i})-f(\eta))\nonumber\\
&\phantom{**}+\lambda_{N}(f(\eta^{N,+})-f(\eta))+\sum_{i=0}^{N}\eta_iq_{i,N}(f(\eta^{i,-})-f(\eta)).
\end{align}

Let $\mu=\bigotimes_{i=1}^{N}\mu_i$ be the product measure of laws $\mu_i$ where for each $i\in \{0,\dots,N\}$, $\mu_i$ is a Poisson law with parameter $C^{(3)}_i$.
We can check that for any configurations  $\eta$, $\eta' \in \Omega$,
\begin{equation}
L(\eta,\eta')\mu(\eta)=L(\eta',\eta)\mu(\eta')
\end{equation}
where $L(\eta,\eta')$ is the transition rate from the configuration $\eta$ to $\eta'$, i.e.,  $L(\eta,\eta')=\mathcal{L}f(\eta)$ with $f(\eta)=\delta_{\eta, \eta'}$.
This implies that the probability measure $\mu$ is reversible and invariant for the Markov process~$\eta$.

Now, consider the model above with injection at rate
$\lambda_{0}$ and absorption at $0$ and only absorption
(without injection) at $N$. 
Such a system can be considered as a $M/G/\infty$ queue where the customers arrive according to a Poisson process of rate $\lambda_0$ and the service time law is that of the lifetime of a particle in the interval $[0,N]$.
Thus, the expected service time of a
customer, denoted by $E[T]$, equals $\Eoiii[\tau_{0}^+\wedge\tau_{N}]$.
By Little's formula (see e.g.\ Section~5.2 of~\cite{Cooper}) we
have
\begin{equation*}
E[T]=\frac{E[R]}{\lambda_{0}}
\end{equation*}
where~$E[R]$ is the mean number of particles in the queue in the
stationary regime.
By a coupling argument, we can see that the distribution of the number of customers in the system in the stationary regime is stochastically dominated by the distribution of the total number of particles in the interval $[0,N]$ in the stationary regime for the particle system with both injection and absorption
of particles at states~$0$ and~$N$. It is not difficult to see that this last distribution is $\mu_0\star\dots\star \mu_N$ (here $\star$ is the convolution product of measures).
Therefore, combining the foregoing observations, we obtain
\begin{equation}
\label{LI3}
\Eoiii[\tau_{0}^+\wedge\tau_{N}] =E[T]=\frac{E[R]}{\lambda_{0}}\leq \frac{1}{\lambda_0}\sum_{x\in \Z} x \mu_0\star\dots\star \mu_N(x) =\frac{1}{C^{(3)}_0}\sum_{x=0}^{N}C^{(3)}_x.
\end{equation}
Finally, by~(\ref{LI1}), (\ref{LI2}) and~(\ref{LI3}) we obtain
\begin{equation*}
\Eo[\tau_B^+\wedge \tau_E] \leq
\frac{1}{C^{(3)}_0\pi_B(0)}\sum_{x=0}^{N}C^{(3)}_x.
\end{equation*}
By Conditions~E and~K, it holds that there exists 
a positive constant~$K_2$ such that $\IP$-a.s.,
\begin{equation*}
\Eo[\tau_B^+\wedge \tau_E]\leq K_2N.
\end{equation*}
This concludes the proof of Lemma~\ref{Little}.
\qed


\section{Proof of Theorem~\ref{Theocond}}
\label{s_proof_Theocond}
In this section, we prove Theorem~\ref{Theocond}. To simplify notations, 
we consider $\sig=1$. Our strategy to prove Theorem \ref{Theocond} is to use 
Theorems~3.6 and~3.10 of~\cite{Dur} 
(which are restated here as Theorems~\ref{Theodu1} and~\ref{Theodu2}). 
These theorems give equivalent conditions for the tightness and convergence of 
finite dimensional distributions of the conditioned processes~$Z^n$ that are easier 
to verify in our case. In~\cite{Dur}, these theorems are stated in a 
quite general form that can be simplified here. 
Also, since in our problem all the processes considered have continuous trajectories, 
we will transpose these theorems on $C[0,1]$ (instead of $D[0,1]$, the Skorokhod space):
\begin{theo}
\label{Theodu1}
The sequence of measures $(\mu_{\omega}^n,n\geq 1)$ is tight if and only if
\begin{equation}
\label{ORTO1}
\lim_{x \to \infty} \limsup_{n\to \infty}\Po[Z^n_1 > x \mid \Lambda_n]=0\phantom{***}\mbox{and}
\end{equation}
\begin{equation}
\label{ORTO2}
\lim_{t \to 0} \limsup_{n\to \infty}\Po[Z^n_t > h \mid \Lambda_n]=0\phantom{**} \mbox{for each $h>0$}.
\end{equation}
\end{theo}
We recall that the measures $\mu_{\omega}^n$ are defined in the introduction.
Now, let us define the following conditions:
\begin{itemize}
\item [(i)] if $x_n\to x$, then $(\Po^{x_n \sqrt{n} }[Z^n_{\cdot}\in \cdot], n\geq 1)$ tends weakly to $P^x[W_{\cdot}\in\cdot]$ in $C[0,1]$,

\item [(ii)] let $x_n\geq 0$, for all $n\geq 1$, then $\lim_{n\to \infty}\Po^{x_n \sqrt{n} }[Z^n_{s}>0, s \leq t_n]=P^x[W_{s}>0, s \leq t]$, whenever $x_n\to x$ and $t_n\to t>0$.

\end{itemize}
\begin{theo}
\label{Theodu2}
Suppose (i)-(ii) hold and $(\mu_{\omega}^n,n\geq 1)$ is tight. Then, $(\mu_{\omega}^n, n\geq 1)$ tends wealky to $W^+$ if and only if 
\begin{equation}
\label{ORTO3}
\lim_{h \to 0} \liminf_{n\to \infty}\Po[Z^n_t > h \mid \Lambda_n]=1\phantom{***}\mbox{for all $t>0$.}
\end{equation}
\end{theo}

In our case, condition (i) is an immediate consequence of the quenched Uniform CLT 
(cf.\ Theorem~1.2 of~\cite{GP}) which in the rest of this paper will be referred as UCLT. 
For condition (ii), let $\eps>0$, we have for all $n$ large enough
\begin{equation*}
\Po^{x_n \sqrt{n} }[Z^n_{s}>0, s \leq t+\eps]\leq \Po^{x_n \sqrt{n} }[Z^n_{s}>0, 
s \leq t_n]\leq \Po^{x_n \sqrt{n} }[Z^n_{s}>0, s \leq t-\eps].
\end{equation*}
Thus, condition (ii) follows from the UCLT and the continuity in $t$ of $P^x[W_{s}>0, 
s \leq t]$.
Our next step is to obtain the weak limit of the sequence $(\Po[Z^n_1 \in \cdot \mid 
\Lambda_n], n\geq1)$. This is the object of Proposition~\ref{propcltord}. 
Then, we obtain the weak limit of $(\Po[Z^n_t \in \cdot \mid \Lambda_n], n\geq 1)$
for all $t\in (0,1)$. 
This is done in Proposition~\ref{prop_Z<x}. 
In the last step, we check that~(\ref{ORTO1}), (\ref{ORTO2}), and~(\ref{ORTO3}) 
hold to end the proof of Theorem~\ref{Theocond}.
\medskip

At this point, 
let us recall some notations of Section~\ref{s_cross_prob_cond}.
Fix $\eps>0$ and define 
$N=\lfloor \eps \sqrt{n} \rfloor$. Let $B=(-\infty,0]$ and
$E=[N,+\infty)$.
Then, define the event $A_{\eps, n}= \{\tau_{E}<\tau^+_B\}$.
For an arbitrary positive integer $M$ define $I_M=[N,N+M]$.
First, let us prove
\begin{prop}
\label{propcltord}
We have $\IP$-a.s.,
\begin{equation}
\label{propordclt}
\lim_{n\to \infty}\Po[Z^n_1 > x \mid \Lambda_n]
=\exp(-x^2/2),\phantom{**}\mbox{for all $x\geq 0$}.
\end{equation}
\end{prop}

\noindent
\textit{Proof.}
For notational convenience,
let us only treat the case $x=1$. 
The generalization to any $x\geq 0$ is straightforward.
Fix $\eps\in (0,1)$, $\delta\in (0,1)$ and 
write
\begin{align}
\label{SUPDEC1}
\Po[X_n >  \sqrt{n}\mid
\Lambda_n]
& =\frac{1}{\Po[\Lambda_n]}\Po[X_n >  \sqrt{n}, A_{\eps,
n},\Lambda_n]\nonumber\\
& =\frac{1}{\Po[\Lambda_n]}\Big(\Po[X_n > \sqrt{n}, A_{\eps,
n},\Lambda_n,X_{\tau_E}\in I_M]
+\Po[X_n >  \sqrt{n}, A_{\eps,
n},\Lambda_n,X_{\tau_E}\notin I_M]\Big)\nonumber\\
&=\frac{1}{\Po[\Lambda_n]}\Big(\Po[X_n >  \sqrt{n}, 
A_{\eps, n},\Lambda_n,X_{\tau_E}\in I_M,\tau_{E}>\delta n]\nonumber\\
& \phantom{*******}+\Po[X_n >  \sqrt{n}, A_{\eps,
n},\Lambda_n,X_{\tau_E}\in I_M,\tau_{E}\leq \delta n]\nonumber\\
& \phantom{*******}+\Po[X_n >  \sqrt{n}, A_{\eps,
n},\Lambda_n,X_{\tau_E}\notin I_M]\Big)\nonumber\\
& =\frac{\Po[A_{\eps, n}]}{\Po[\Lambda_n]}\Big(\Po[  X_{\tau_E}\in
I_M\mid A_{\eps, n}] \Po[ \tau_{E}>\delta n \mid X_{\tau_E}\in
I_M,A_{\eps, n}]\nonumber\\
& \phantom{********}\times \Po[X_n > \sqrt{n},\Lambda_n \mid
X_{\tau_E}\in I_M, A_{\eps, n},\tau_{E}>\delta n]  \nonumber\\
& \phantom{*******}+\Po[X_{\tau_E}\in I_M \mid A_{\eps, n}] \Po[
\tau_{E}\leq \delta n \mid X_{\tau_E}\in I_M, A_{\eps, n}]
\nonumber\\
& \phantom{********}\times \Po[X_n > \sqrt{n},\Lambda_n \mid
X_{\tau_E}\in I_M, A_{\eps, n},\tau_{E}\leq \delta n]\nonumber\\
& \phantom{*******}+\Po[X_n > 
\sqrt{n},\Lambda_n,X_{\tau_E}\notin I_M\mid A_{\eps,n}]\Big).
\end{align}
Informally, the rest of the proof consists in 
using the decomposition~(\ref{SUPDEC1}) in order to find good lower
and upper bounds~$L_n$ and~$U_n$ for $\Po[X_n >  \sqrt{n}\mid
\Lambda_n]$ such that $U_n/L_n\to 1$ as $n \to \infty$.
We start with the upper bound.
Let us write 
\begin{align}
\label{UPCLT}
\Po[X_n >  \sqrt{n}\mid \Lambda_n]
& \leq \frac{\Po[A_{\eps, n}]}{\Po[\Lambda_n]}\Big(\Po[X_{\tau_E}\notin I_M\mid A_{\eps, n}] +
\Po[\tau_{E}>\delta n \mid X_{\tau_E}\in I_M, A_{\eps, n}]\nonumber\\
& \phantom{*********}+\Po[X_n >  \sqrt{n},\Lambda_n \mid X_{\tau_E}\in
I_M, A_{\eps, n},\tau_{E}\leq \delta n] 
\Big).
\end{align}
Observe that we can bound the term $\Po[X_{\tau_E}\notin I_M\mid
A_{\eps, n}]$
from above using Lemma~\ref{SNLcond}: let $\eta>0$, then we can choose $M$ large enough in such a way that
\begin{equation}
\label{UPBDRT1}
\Po[X_{\tau_E}\notin I_M\mid A_{\eps,n}] \leq \eta.
\end{equation}
 Next, let us bound the other terms of the right-hand side
of~(\ref{UPCLT}) from above.
For $\Po[A_{\eps,n}]/\Po[\Lambda_n]$, we write
\begin{align}
\label{Lowlamb}
\Po[\Lambda_n]
 &\geq \Po[\Lambda_n,A_{\eps,n},X_{\tau_E}\in I_M,\tau_{E}\leq \delta
n]\nonumber\\
 &= \Po[A_{\eps,n}]\Po[X_{\tau_E}\in I_M\mid A_{\eps,n}]
\Po[\tau_{E}\leq \delta n \mid X_{\tau_E}\in I_M, A_{\eps,n}]
\Po[\Lambda_n\mid A_{\eps,n},X_{\tau_E}\in
I_M,\tau_{E}\leq \delta n].
\end{align}
Hence,
\begin{align*}
\frac{\Po[\Lambda_n]}{\Po[A_{\eps,n}]}
 &\geq \Po[X_{\tau_E}\in I_M\mid A_{\eps,n}] \Po[\tau_{E}\leq \delta
n \mid X_{\tau_E}\in I_M, A_{\eps,n}] 
\Po[\Lambda_n\mid A_{\eps,n},X_{\tau_E}\in
I_M,\tau_{E}\leq \delta n].
\end{align*}
Again, we use Lemma~\ref{SNLcond} to bound the term
$\Po[X_{\tau_E}\in I_M\mid A_{\eps,n}] $ from below.
 For the term $\Po[\tau_{E}\leq \delta n\mid X_{\tau_E}\in
I_M,A_{\eps,n}]$ we write
\begin{align}
\label{REL1}
\Po[\tau_{E}\leq \delta n\mid X_{\tau_E}\in
I_M,A_{\eps,n}]
 &=1-\Po[\tau_{E}> \delta n\mid X_{\tau_E}\in I_M,
A_{\eps,n}]
\end{align}
and
\begin{align}
\label{EQ1}
\Po[ \tau_{E}>\delta n \mid X_{\tau_E}\in I_M,A_{\eps, n}]
 &=\frac{\Po[ \tau_{E}>\delta n, X_{\tau_E}\in I_M, A_{\eps,
n}]}{\Po[X_{\tau_E}\in I_M, A_{\eps, n}]} =\frac{\Po[ \tau_{E}>\delta n, X_{\tau_E}\in I_M, A_{\eps,
n}]}{\Po[X_{\tau_E}\in I_M\mid A_{\eps, n}] \Po[A_{\eps, n}]}.
\end{align}
We first treat the numerator of~(\ref{EQ1}). By Chebyshev's
inequality we obtain
\begin{equation*}
 \Po[ \tau_{E}>\delta n,X_{\tau_E}\in I_M, A_{\eps, n}]\leq
\Po[\tau_B^+\wedge \tau_{E}>\delta n]\leq \frac{\Eo[\tau_B^+\wedge
\tau_{E}]}{\delta n}.
\end{equation*}
Using~(\ref{EQ1}) and Lemmas~\ref{Little}, \ref{SNLcond} 
and~\ref{Aen} we obtain 
\begin{equation}
\label{FeD1}
 \Po[ \tau_{E}>\delta n \mid X_{\tau_E}\in I_M, A_{\eps, n}]\leq
\frac{K_2N^2}{K_1\delta n (1-\eta)}.
\end{equation}
Then, we deal with the term $\Po[\Lambda_n\mid
A_{\eps,n},X_{\tau_E}\in
I_M,\tau_{E}\leq \delta n]$.
By the Markov property we obtain
\begin{align}
\label{Markuv}
\lefteqn{
\Po[\Lambda_n\mid A_{\eps,n},X_{\tau_E}\in
I_M,\tau_{E}\leq
\delta n]
}\nonumber\\
 &= \frac{1}{\Po[A_{\eps, n},X_{\tau_E}\in I_M,\tau_{E}\leq \delta
n]} \sum_{x\in I_M} \sum_{u=1}^{\lfloor \delta n \rfloor}
\Po[\Lambda_n \mid X_{\tau_E}=x,\tau_{E}=u, A_{\eps, n}]
\Po[X_{\tau_E}=x,\tau_{E}=u,A_{\eps, n}]\nonumber\\ 
 &\geq \min_{x\in I_M}\min_{u\leq \lfloor \delta n
\rfloor}\Po^{x}[\Lambda_{n-u}]\nonumber\\
&\geq \min_{x\in I_M} \Po^{x}[\Lambda_{n}].
\end{align}
Thus, by~(\ref{Lowlamb}), (\ref{REL1}),  (\ref{EQ1}), (\ref{Markuv}) and 
Lemma~\ref{SNLcond}, we have
\begin{equation}
\label{UPBDRT2}
 \frac{\Po[\Lambda_n]}{\Po[A_{\eps,n}]}\geq
(1-\eta)\Big(1-\frac{K_2N^2}{K_1\delta n (1-\eta)}\Big)\min_{x\in
I_M} \Po^{x}[\Lambda_{n}].
\end{equation}
 To bound the term $\Po[X_n >\sqrt{n},\Lambda_n \mid
X_{\tau_E}\in I_M, A_{\eps, n},\tau_{E}\leq \delta n]$ from above we
do the following. 
Let us denote by
$\mathcal{E}$ the event 
$\{X_{\tau_E}\in I_M, A_{\eps, n},\tau_{E}\leq \delta n\}$.
Since $A_{\eps, n} \in {\mathcal F}_{\tau_{E}}$ the $\sigma$-field
generated by~$X$ until the stopping time~$\tau_E$, we have by the
Markov property and the fact that $\delta<1$,
\begin{align}
\label{TYTY}
\Po[X_n>  \sqrt{n},\Lambda_n \mid
\mathcal{E}]
 &= \frac{1}{\Po[\mathcal{E}]}  \sum_{x\in I_M}\sum_{u=1}^{\lfloor
\delta n \rfloor} \Po[X_n >  \sqrt{n},\Lambda_n \mid
X_{\tau_E}=x,\tau_{E}=u, A_{\eps, n}] 
\nonumber\\
 &\phantom{***********}\times 
\Po[X_{\tau_E}=x,\tau_{E}=u,A_{\eps,
n}]\nonumber\\
&\leq \max_{x\in I_M} \max_{u \leq \lfloor \delta n \rfloor}
\Po[X_n >  \sqrt{n},\Lambda_n \mid
X_{\tau_E}=x,\tau_{E}=u,A_{\eps, n}]\nonumber\\
 &=  \max_{x\in I_M} \max_{u \leq \lfloor \delta n
\rfloor}\Po^{x}[X_{n-u}>  \sqrt{n},X_k>0, 1\leq k\leq
n-u]\nonumber\\
 &= \max_{x\in I_M} \max_{u \leq \lfloor \delta n
\rfloor}\Po^{x}[X_{n-u} > \sqrt{n},\Lambda_{n-u}].
\end{align}
Now, fix $\delta'\in (0,1)$. 
Then, we use the following estimate for $x\in I_M$ and $u \leq
\lfloor \delta n \rfloor$,
\begin{align*}
\Po^{x}[X_{n-u} > 
\sqrt{n},\Lambda_{n-u}]
 &\leq \Po^{x}\big[(\{X_{n-\lfloor \delta n \rfloor} > (1-\delta')
\sqrt{n}\}\cup\{|X_{n-\lfloor \delta n \rfloor}-X_{n-u}|> \delta' 
\sqrt{n}\})\cap\Lambda_{n-u}\big]\nonumber\\
 &\leq  \Po^{x}\Big[\Big(\Big\{X_{n-\lfloor \delta n \rfloor} >
(1-\delta') \sqrt{n}\Big\}\nonumber\\
 &\phantom{*****}\cup\Big\{\max_{u \leq \lfloor \delta n
\rfloor}|X_{n-\lfloor \delta n \rfloor}-X_{n-u}|> \delta' 
\sqrt{n}\Big\}\Big)\cap\Lambda_{n-\lfloor \delta n
\rfloor}\Big].\nonumber\\
\end{align*}
Hence, we obtain for all $x\in I_M$ that
\begin{align}
\label{DEC1}
\max_{u \leq \lfloor \delta n \rfloor}\Po^{x}[X_{n-u} >
 \sqrt{n},\Lambda_{n-u}]
 &\leq \Po^{x}[X_{n-\lfloor \delta n \rfloor} > (1-\delta')
\sqrt{n},\Lambda_{n-\lfloor \delta n \rfloor}]\nonumber\\
 &\phantom{**}+\Po^{x}\Big[\max_{u \leq \lfloor \delta n
\rfloor}|X_{n-\lfloor \delta n \rfloor}-X_{n-u}|> \delta' 
\sqrt{n},\Lambda_{n-\lfloor \delta n \rfloor}\Big].
\end{align}
To sum up, using~(\ref{UPBDRT1}), (\ref{UPBDRT2}), (\ref{FeD1}) 
and~(\ref{DEC1}) we obtain that $\IP$-a.s., 
\begin{align}
\label{UPBOUNDORD}
\Po[X_n>\sqrt{n}\mid \Lambda_n]
 &\leq (1-\eta)^{-1}\Big(1-\frac{K_2N^2}{K_1\delta n
(1-\eta)}\Big)^{-1}\Big(\min_{x\in I_M}
\Po^{x}[\Lambda_{n}]\Big)^{-1}\Big(\frac{K_2N^2}{K_1\delta n
(1-\eta)}+\eta \nonumber\\
 &\phantom{**}+\max_{x\in I_M}\Po^{x}[X_{n-\lfloor \delta n \rfloor}
> (1-\delta') \sqrt{n},\Lambda_{n-\lfloor \delta n
\rfloor}]\nonumber\\
 &\phantom{**}+\max_{x\in I_M} \Po^{x}\Big[\max_{u \leq \lfloor
\delta n \rfloor}|X_{n-\lfloor \delta n \rfloor}-X_{n-u}|> \delta' 
\sqrt{n},\Lambda_{n-\lfloor \delta n \rfloor}\Big]\Big).
\end{align}

 Our goal is now to calculate the $\limsup$ as $n \to \infty$ of
both sides of~(\ref{UPBOUNDORD}).
Let us first compute $\limsup_{n\to
\infty}(\Po^{x}[\Lambda_{n}])^{-1}$ for $x\in I_M$. We have by
definition of $Z^n$
\begin{equation*}
 \Po^{x}[\Lambda_{n}]=\Po^{x}[X_m>0,0\leq m \leq
n]=\Po^{x}\big[Z^n_t>0, t\in [0,1]\big].
\end{equation*}
Thus, by the UCLT, we have 
\begin{equation*}
 \lim_{n\to \infty}\Po^{x}\big[Z^n_t>0, t\in
[0,1]\big]=P^{\eps}\Big[\min_{0\leq t\leq 1}W(t)>0\Big]
\end{equation*}
with~$W$ a standard Brownian motion.
Using the reflexion principle (see Chap.\ III, Prop.\ 3.7 in \cite{RY}), we obtain
\begin{align*}
P^{\eps}\Big[\min_{0\leq t\leq 1}W(t)>0\Big]
&=P^0[|W(1)|< \eps]
=\int_{-\eps}^{\eps}\frac{1}{\sqrt{2\pi}}
e^{-\frac{x^2}{2}}dx.
\end{align*}
So, we obtain
\begin{equation}
\label{URT12}
 \lim_{n\to \infty}\min_{x\in I_M} (\Po^{x}[\Lambda_{n}])^{-1}=
\Big(\int_{-\eps}^{\eps}\frac{1}{\sqrt{2\pi}}e^{-
\frac{x^2}{2}}dx\Big)^{-1}=\Big(\frac{2\eps}{\sqrt{2\pi}}
+o(\eps)\Big)^{-1}
\end{equation}
as $\eps\to 0$.

Now, let us bound 
$\limsup_{n \to \infty} \Po^{x}[X_{n-\lfloor \delta n \rfloor} >
(1-\delta') \sqrt{n},\Lambda_{n-\lfloor \delta n \rfloor}]$ from
above. We have
\begin{align*}
\Po^{x}[X_{n-\lfloor \delta n \rfloor} > (1-\delta')
\sqrt{n},\Lambda_{n-\lfloor \delta n
\rfloor}]
 &\leq \Po^{x}\Big[X_{n-\lfloor \delta n \rfloor} > (1-\delta')
\sqrt{n-\lfloor \delta n \rfloor},\Lambda_{n-\lfloor \delta n
\rfloor}\Big]\nonumber\\
 &=\Po^{x}\Big[Z^{n-\lfloor \delta n \rfloor}_1 >
(1-\delta'),Z^{n-\lfloor \delta n \rfloor}_t>0, t\in
[0,1]\Big].
\end{align*}
 As $\delta<1$ and $x\in I_M$, we have by the UCLT,
\begin{align*}
 \lim_{n\to \infty}\Po^{x} \Big[Z^{n-\lfloor \delta n \rfloor}_1 > (1-\delta')
,Z^{n-\lfloor \delta n \rfloor}_t>0, t\in
[0,1]\Big]=
 &P^{\frac{\eps}{\sqrt{1-\delta}}}\Big[W(1) > (1-\delta')
,\min_{0\leq t\leq 1}W(t)>0\Big].
\end{align*}
 Abbreviate $\eps':=\eps (1-\delta)^{-\frac{1}{2}}$ and let us compute
$P^{\eps'}\Big[W(1) > (1-\delta'),\min_{0\leq t\leq
1}W(t)>0\Big]$ for sufficiently small~$\eps$.
By the reflexion principle for Brownian motion, we have
\begin{align*}
P^{\eps'}\Big[W(1) > (1-\delta')
,\min_{0\leq t\leq 1}W(t)>0\Big]
 &=P^{\eps'}\Big[W(1) >( 1-\delta')
\Big]-P^{\eps'}\Big[W(1) <- (1-\delta')\Big]\nonumber\\
 &=P\Big[W(1)> 1-(\delta'+\eps') \Big]-P\Big[W(1)
<-1+(\delta'-\eps') )\Big]\\
&=\frac{1}{\sqrt{2\pi}}\int_{1-(\delta'+\eps')}^
{1-(\delta'-\eps')}e^{-\frac{x^2}{2}}dx.
\end{align*}
Therefore, we obtain, as $\eps\to 0$
\begin{align}
\label{SpL2}
 \limsup_{n \to \infty} \max_{x\in I_M} \Po^{x}[X_{n-\lfloor \delta n
\rfloor} > (1-\delta') \sqrt{n},\Lambda_{n-\lfloor \delta n
\rfloor}]
 &\leq \frac{1}{\sqrt{2\pi}}\int_{1-(\delta'+\eps')}^{
1-(\delta'-\eps')}e^{-\frac{x^2}{2}}dx\nonumber\\
 &=\frac{2\eps}{\sqrt{2\pi(1-\delta)}}e^{-\frac{1}{2}}+o(\eps).
\end{align}

 Then, let us bound $\limsup_{n \to \infty} \Po^{x}\Big[\max_{u \leq
\lfloor \delta n \rfloor}|X_{n-\lfloor \delta n \rfloor}-X_{n-u}|>
\delta'  \sqrt{n},\Lambda_{n-\lfloor \delta n \rfloor}\Big]$ from
above in (\ref{UPBOUNDORD}) for $x\in I_M$. First, observe that
\begin{align*}
\Po^{x}\Big[\max_{u \leq \lfloor \delta n
\rfloor}|X_{n-\lfloor \delta n \rfloor}-X_{n-u}|\geq \delta' 
\sqrt{n},\Lambda_{n-\lfloor \delta n
\rfloor}\Big]
 &\leq \Po^{x}\Big[\max_{u \leq \lfloor \delta n
\rfloor}|X_{n-\lfloor \delta n \rfloor}-X_{n-u}|\geq \delta' 
\sqrt{n}\Big]
\end{align*}
and
\begin{align*}
\lefteqn{
\Po^{x}\Big[\max_{u \leq \lfloor \delta n
\rfloor}|X_{n-\lfloor \delta n \rfloor}-X_{n-u}|\geq \delta' 
\sqrt{n}\Big]
}
\nonumber\\
 &=\Po^{x}\Big[\max_{n-\lfloor \delta n \rfloor\leq k\leq
n}|X_k-X_{n-\lfloor \delta n \rfloor}|\geq \delta' 
\sqrt{n}\Big]\nonumber\\
 &\leq \Po^{x}\Big[\max_{n-\lfloor \delta n \rfloor\leq k\leq
n}(X_k-\min_{n-\lfloor \delta n \rfloor\leq l\leq k}X_l) \geq \delta'
 \sqrt{n}\Big]
+\Po^{x}\Big[\min_{n-\lfloor \delta n \rfloor\leq k\leq
n}(X_k-\max_{n-\lfloor \delta n \rfloor\leq l\leq k}X_l) \leq-\delta'
 \sqrt{n}\Big]\nonumber\\
 &\leq \Po^{x}\Big[\max_{1-\delta \leq t\leq
1}(Z_t^n-\min_{1-\delta\leq s\leq t}Z_s^n) \geq \delta' 
\Big]
+\Po^{x}\Big[\min_{1-\delta \leq
t\leq 1}(Z_t^n-\max_{1-\delta\leq s\leq t}Z_s^n) \leq-\delta' 
\Big].
\end{align*}
Using the UCLT, we
obtain
\begin{align}
\label{Convfraca1}
\lim_{n\to \infty}\Po^{x}\Big[\max_{1-\delta \leq t\leq 1}(Z_t^n-\min_{1-\delta\leq
s\leq t}Z_s^n) \geq \delta' \Big]=
 &P^{\eps}\Big[\max_{1-\delta \leq t\leq
1}\big(W(t)-\min_{1-\delta\leq s\leq t}W(s)\big) \geq \delta' \Big]
\end{align}
and 
\begin{align}
\label{Convfraca2}
\lim_{n\to \infty} \Po^{x}\Big[\min_{1-\delta \leq t\leq 1}(Z_t^n-\max_{1-\delta\leq
s\leq t}Z_s^n) \leq-\delta' \Big]=
 &P^{\eps}\Big[\min_{1-\delta \leq t\leq
1}\big(W(t)-\max_{1-\delta\leq s\leq t}W(s)\big) \leq-\delta' \Big].
\end{align}
 Observe that the right-hand sides of~(\ref{Convfraca1}) 
and~(\ref{Convfraca2}) are equal since $(-W)$ is a Brownian motion.
Thus, let us compute for example $P^{\eps}[\max_{1-\delta
\leq t\leq 1}(W(t)-\min_{1-\delta\leq s\leq t}W(s)) \geq \delta' ]$.
 First, by the Markov property and 
since the event is invariant by space shifts, we have 
\begin{align*}
P^{\eps}\Big[\max_{1-\delta \leq t\leq
1}\big(W(t)-\min_{1-\delta\leq s\leq t}W(s)\big) \geq \delta' \Big]
 &=P\Big[\max_{0 \leq t\leq \delta}\big(W(t)-\min_{0 \leq s\leq t}W(s)\big)
\geq \delta' \Big].
\end{align*}
 By L\'evy's Theorem (cf.~\cite{RY}, Chapter~VI, Theorem~2.3), we
have
\begin{equation*}
 P\Big[\max_{0 \leq t\leq \delta}\big(W(t)-\min_{0 \leq s\leq t}W(s)\big)
\geq \delta' \Big]=P\Big[\max_{0\leq t\leq
\delta}|W(t)|\geq \delta' \Big].
\end{equation*}
Then, by the reflexion principle, we have
\begin{equation*}
P\Big[\max_{0\leq t\leq \delta}|W(t)|\geq \delta' \Big]
 \leq 2P\Big[\max_{0\leq t\leq \delta}W(t)\geq \delta'\Big]
=4P[W(\delta)\geq \delta' ].
\end{equation*}
 Using an estimate on the tail of the Gaussian law
(cf.~\cite{PerMot}, Appendix~II, Lemma~3.1) we obtain
\begin{equation*}
 P\Big[\max_{0\leq t\leq \delta}|W(t)|\geq \delta' \Big]\leq
\frac{4\sqrt{\delta}}{\delta'
\sqrt{2\pi}}\exp\Big\{-\frac{(\delta')^2}{2\delta}\Big\}.
\end{equation*}
Thus, we find
\begin{align}
\label{SpiL}
\limsup_{n\to \infty} \max_{x\in I_M}\Po^{x}\Big[\max_{u
\leq \lfloor \delta n \rfloor}|X_{n-\lfloor \delta n
\rfloor}-X_{n-u}|> \delta'  \sqrt{n},\Lambda_{n-\lfloor \delta n
\rfloor}\Big]
 &\leq \frac{8\sqrt{\delta}}{\delta'
\sqrt{2\pi}}\exp\Big\{-\frac{(\delta')^2}{2\delta}\Big\}.
\end{align}
Finally, combining~(\ref{UPBOUNDORD}), (\ref{URT12}), (\ref{SpL2})
and~(\ref{SpiL}), we obtain
\begin{align}
\label{UPBOUNDORD2}
\limsup_{n \to \infty}\Po[X_n> \sqrt{n}\mid
\Lambda_n]
 &\leq (1-\eta)^{-1}\Big(1-\frac{K_2\eps^2}{K_1\delta
(1-\eta)}\Big)^{-1}\Big(\frac{2\eps}{\sqrt{2\pi}}+o(\eps)\Big)^{
-1}\Big(\frac{K_2\eps^2}{K_1\delta  (1-\eta)}+\eta \nonumber\\
&\phantom{*******}+ \frac{2\eps}{\sqrt{2\pi(1-\delta)}}
e^{-\frac{1}{2}}+o(\eps)
+ \frac{8\sqrt{\delta}}{\delta' \sqrt{2\pi}}
\exp\Big\{-\frac{(\delta')^2}{2\delta}\Big\}                     
\Big).
\end{align}

 Next, let us bound the quantity $\Po[X_n > \sqrt{n}\mid
\Lambda_n]$ from below. Using (\ref{SUPDEC1}), we write 
\begin{align}
\label{RUIT0}
\Po[X_n > \sqrt{n}\mid \Lambda_n]
 &\geq \frac{\Po[A_{\eps, n}]}{\Po[\Lambda_n]}\Po[X_{\tau_E}\in
I_M\mid A_{\eps, n}]\Po[ \tau_{E}\leq\delta n \mid X_{\tau_E}\in I_M,
A_{\eps, n}] \nonumber\\
 &\phantom{**}\times \Po[X_n >  \sqrt{n},\Lambda_n \mid
X_{\tau_E}\in I_M, A_{\eps, n},\tau_{E}\leq \delta n].
\end{align}
 As we have already treated the terms $\Po[ \tau_{E}\leq\delta n \mid
X_{\tau_E}\in I_M, A_{\eps, n}]$ and $\Po[X_{\tau_E}\in
I_M\mid A_{\eps, n}]$ in~(\ref{REL1}) and Lemma \ref{SNLcond} respectively, we just need to
bound the terms  $\Po[A_{\eps, n}]/\Po[\Lambda_n]$ and 
$\Po[X_n >  \sqrt{n},\Lambda_n \mid X_{\tau_E}\in I_M, A_{\eps, n},\tau_{E}\leq
\delta n]$ from below.

Let us start with the term $\Po[A_{\eps, n}]/\Po[\Lambda_n]$. Observe
that
\begin{align}
\label{Mo0}
\Po[\Lambda_n]
 &=\Po[\Lambda_n,\tau_{E}\leq \delta n]+\Po[\Lambda_n,\tau_{E}>
\delta n]\nonumber\\
  &=\Po[\Lambda_n,A_{\eps,n},\tau_{E}\leq \delta
n]+\Po[\Lambda_n,A_{\eps,n},\tau_{E}> \delta n]
+\Po[\Lambda_n,A^c_{\eps,n},\tau_{E}> \delta
n]\nonumber\\
 &\leq \Po[\Lambda_n,A_{\eps,n},\tau_{E}\leq \delta
n]+\Po[\Lambda_n,A_{\eps,n},\tau_{E}> \delta
n]+\Po[\Lambda_n,A^c_{\eps,n}]\nonumber\\
 &\leq \Po[\Lambda_n,A_{\eps,n},\tau_{E}\leq \delta n,X_{\tau_E}\in
I_M]+\Po[\Lambda_n,A_{\eps,n},\tau_{E}> \delta n,X_{\tau_E}\in
I_M]\nonumber\\
 &\phantom{**}+2 \Po[X_{\tau_E}\notin I_M,
A_{\eps,n}]+\Po[\Lambda_n,A^c_{\eps,n}]\nonumber\\
 &\leq \Po[A_{\eps,n}] \Big[\Po[\Lambda_n\mid A_{\eps,n},\tau_{E}\leq
\delta n,X_{\tau_E}\in I_M]+2 \Po[X_{\tau_E}\notin I_M\mid
A_{\eps,n}]\nonumber\\
 &\phantom{**}+\Po[\tau_{E}> \delta n\mid X_{\tau_E}\in I_M,
A_{\eps,n}]+\frac{\Po[\Lambda_n,A^c_{\eps,n}]}{\Po[A_{\eps,n}]}\Big].
\end{align}
 From the first equality in~(\ref{Markuv}) we obtain
\begin{equation}
\label{MO1}
 \Po[\Lambda_n\mid A_{\eps,n},X_{\tau_E}\in I_M,\tau_{E}\leq \delta
n]
\leq \max_{x\in I_M}\max_{u\leq \lfloor \delta n
\rfloor}\Po^{x}[\Lambda_{n-u}]
\leq \max_{x\in I_M} \Po^{x}[\Lambda_{n-\lfloor \delta n \rfloor}].
\end{equation}
Now, let us treat the term $\Po[\Lambda_n, A^c_{\eps,n}]$.
First, observe that by definition of $A_{\eps,n}$ we have
\[\Po[\Lambda_n, A^c_{\eps,n}]\leq \Po[\tau_B^+\wedge \tau_E>n].
\]
Then, by Chebyshev's inequality we obtain
\begin{equation*}
 \Po[\tau_B^+\wedge \tau_E>n]\leq \frac{\Eo[\tau_B^+\wedge
\tau_E]}{n}.
\end{equation*}
By Lemma~\ref{Little}, we obtain
\begin{equation}
\label{MO2}
\Po[\tau_B^+\wedge \tau_E>n]\leq \frac{K_2N}{n}.
\end{equation}
Thus, by~(\ref{FeD1}), (\ref{Mo0}), (\ref{MO1}), (\ref{MO2}) and
Lemmas~\ref{SNLcond}
and~\ref{Aen} we obtain
\begin{align}
\label{RUIT1}
\frac{\Po[A_{\eps, n}]}{\Po[\Lambda_n]}
& \geq \Big(\max_{x\in I_M} \Po^x[\Lambda_{n-\lfloor \delta n
\rfloor}]+2\eta + \frac{K_2N^2}{K_1\delta n (1-\eta)}+\frac{K_2
N^2}{K_1n}\Big)^{-1}.
\end{align}

 Let us find a lower bound for $\Po[X_n > \sqrt{n},\Lambda_n
\mid X_{\tau_E}\in I_M, A_{\eps, n},\tau_{E}\leq \delta n]$ in (\ref{RUIT0}).
 Since $A_{\eps, n} \in {\mathcal F}_{\tau_{E}}$ we have by the
Markov property,
\begin{align}
 \lefteqn{\Po[X_n >  \sqrt{n},\Lambda_n \mid X_{\tau_E}\in
I_M,A_{\eps, n},\tau_{E}\leq \delta
n]}\phantom{*************}\nonumber\\
 &\geq \min_{x\in I_M}\min_{u \leq \lfloor \delta n \rfloor}\Po[X_n >
 \sqrt{n},\Lambda_n \mid X_{\tau_E}=x,\tau_{E}=u,A_{\eps,
n}]\nonumber\\
 &= \min_{x\in I_M}\min_{u \leq \lfloor \delta n
\rfloor}\Po^{x}[X_{n-u} >  \sqrt{n},X_k>0, 1\leq k\leq
n-u]\nonumber\\
 &= \min_{x\in I_M}\min_{u \leq \lfloor \delta n
\rfloor}\Po^{x}[X_{n-u} >  \sqrt{n},\Lambda_{n-u}].
\end{align}
 For $x\in I_M$ and 
$u \leq \lfloor \delta n \rfloor$ we 
write 
\begin{align}
\label{RUIT2}
\Po^{x}[X_{n-u} >  \sqrt{n},\Lambda_{n-u}]
 &\geq \Po^{x}[X_{n} > (1+\delta') \sqrt{n},|X_n-X_{n-u}|\leq
\delta'\sqrt{n},\Lambda_{n-u}]\nonumber\\
 &\geq \Po^{x}\Big[X_{n} > (1+\delta') \sqrt{n},\max_{u \leq
\lfloor \delta n \rfloor}|X_n-X_{n-u}|\leq
\delta'\sqrt{n},\Lambda_{n-u}\Big]\nonumber\\
  &\geq \Po^{x}\Big[X_{n} > (1+\delta') \sqrt{n},\max_{u \leq
\lfloor \delta n \rfloor}|X_n-X_{n-u}|\leq
\delta'\sqrt{n},\Lambda_{n}\Big]\nonumber\\
  &\geq \Po^{x}[X_{n} > (1+\delta')
\sqrt{n},\Lambda_{n}]-\Po^{x}\Big[\max_{u \leq \lfloor \delta n
\rfloor}|X_n-X_{n-u}|> \delta'\sqrt{n}\Big].
\end{align}
To sum up, by~(\ref{RUIT0}), (\ref{RUIT1}), (\ref{RUIT2}),
(\ref{REL1}), (\ref{FeD1}) and Lemma~\ref{SNLcond} we obtain that $\IP$-a.s., 
\begin{align}
\label{FGH}
\Po[X_n> \sqrt{n}\mid \Lambda_n]
 &\geq (1-\eta)\Big(1-\frac{K_2N^2}{K_1\delta n
(1-\eta)}\Big)\nonumber\\
 &\phantom{**}\times \Big(\max_{x\in I_M} \Po^x[\Lambda_{n-\lfloor
\delta n
\rfloor}]+2\eta + \frac{K_2N^2}{K_1\delta n (1-\eta)}+\frac{K_2
N^2}{K_1n}\Big)^{-1}\nonumber\\
 &\phantom{**}\times \Big(\min_{x\in I_M}\Po^{x}[X_{n} >
(1+\delta')
\sqrt{n},\Lambda_{n}]-\max_{x\in I_M} \Po^{x}\Big[\max_{u \leq
\lfloor \delta n \rfloor}|X_n-X_{n-u}|> \delta'\sqrt{n}\Big]\Big).
\end{align}

 Let us now compute $\liminf_{n\to \infty}$ of both sides
of~(\ref{FGH}). 
First, by~(\ref{URT12}) we have
\begin{equation}
\label{URT26}
 \lim_{n\to \infty}\max_{x\in I_M} \Po^{x}[\Lambda_{n-\lfloor \delta
n \rfloor}]=\frac{2\eps}{\sqrt{2\pi(1-\delta)}}+o(\eps)
\end{equation}
as $\eps\to 0$.
Then, by the UCLT and after some elementary 
computations similar to those which led to (\ref{SpL2}) and (\ref{SpiL}) we obtain
\begin{equation}
\label{UTI1}
\lim_{n \to \infty} \min_{x\in I_M} \Po^{x}[X_{n} >
(1+\delta')
\sqrt{n},\Lambda_{n}]
=\frac{1}{\sqrt{2\pi}}\int_{1+(\delta'-\eps)}^{
1+(\delta'+\eps)}e^{-\frac{x^2}{2}}dx
=\frac{2\eps}{\sqrt{2\pi}}e^{-\frac{1}{2}}+o(\eps)
\end{equation}
as $\eps\to 0$, and 
\begin{align}
\label{UTI2}
\limsup_{n \to \infty}\max_{x\in I_M} \Po^{x}\Big[\max_{u
\leq \lfloor \delta n \rfloor}|X_n-X_{n-u}|>
\delta'\sqrt{n}\Big]
 &\leq \frac{8\sqrt{\delta}}{\delta'
\sqrt{2\pi}}\exp\Big\{-\frac{(\delta')^2}{2\delta}\Big\}.
\end{align}
Thus, combining~(\ref{FGH}) with (\ref{URT26}), (\ref{UTI1}) 
and~(\ref{UTI2}) leads to 
\begin{align}
\label{UPBOUNDORD3}
\liminf_{n \to \infty}\Po[X_n> \sqrt{n}\mid \Lambda_n]
 &\geq (1-\eta)\Big(1-\frac{K_2\eps^2}{K_1\delta 
(1-\eta)}\Big)\nonumber\\
 &\phantom{**}\times \Big(
\frac{2\eps}{\sqrt{2\pi(1-\delta)}}+o(\eps)+2\eta +
\frac{K_2\eps^2}{K_1\delta (1-\eta)}+\frac{K_2
\eps^2}{K_1}\Big)^{-1}\nonumber\\
 &\phantom{**}\times
\Big(\frac{2\eps}{\sqrt{2\pi}}e^{-\frac{1}{2}}+o(\eps)-
\frac{8\sqrt{\delta}}{\delta'\sqrt{2\pi}}
\exp\Big\{-\frac{(\delta')^2}{2\delta}\Big\} \Big).
\end{align}

 Now take $\eta= \eps^2$, $\delta=\eps^{\frac{1}{2}}$ and
$\delta'=\eps^{\frac{1}{8}}$ and let $\eps \to 0$ in~(\ref{UPBOUNDORD2}) 
and~(\ref{UPBOUNDORD3}) to prove~(\ref{propordclt}).
\qed

 The next step is to show the weak convergence of $(\Po[Z^n_t \in \cdot \mid \Lambda_n], n\geq 1)$ for all $t\in (0,1)$. 
We start by recalling the transition density function from $(0,0)$ to $(t,y)$ of the Brownian meander (see \cite{Igle}):
 \begin{equation}
\label{GRANFIN}
q(t,y)=t^{-\frac{3}{2}}y
\exp\Big(-\frac{y^2}{2t}\Big){\tilde N}(y(1-t)^{-\frac{1}{2}})
\end{equation}
for $y>0$, $0<t \leq 1$,
where
\[
 {\tilde N}(x)=\sqrt{\frac{2}{\pi}}\int_0^x
e^{-\frac{u^2}{2}}du
\]
for $x\geq 0$.
We will prove the following
\begin{prop}
\label{prop_Z<x}
We have $\IP$-a.s., for all $x\geq 0$ and $0<t<1$,
\begin{equation}
\label{FFD1}
 \lim_{n\to \infty}\Po[Z^n_t\leq x\mid \Lambda_n]=\int_0^x
q(t,y)dy.
\end{equation}
\end{prop}
\textit{Proof.}
First notice the following. For all $\tilde{\eps}>0$ we have
\begin{align}
\label{Finiteconvdist}
\lefteqn{\Po\Big[Z^n_{\frac{\lfloor nt\rfloor}{n}}
\leq x-\tilde{\eps} \mid \Lambda_n\Big]}\phantom{******}\nonumber\\
&\leq \Po\Big[ Z^n_{\frac{\lfloor nt\rfloor}{n}}\leq x-\tilde{\eps}, \Big| 
Z^n_{\frac{\lfloor nt\rfloor+1}{n}}-Z^n_{\frac{\lfloor
nt\rfloor}{n}}\Big|\leq \tilde{\eps} \mid \Lambda_n\Big]+\Po\Big[\Big|
Z^n_{\frac{\lfloor nt\rfloor+1}{n}}-Z^n_{\frac{\lfloor
nt\rfloor}{n}}\Big|> \tilde{\eps} \mid \Lambda_n\Big]\nonumber\\
&\leq \Po[Z^n_t\leq x\mid \Lambda_n]
+\Po[\Lambda_n]^{-1}\Po[|X_{\lfloor nt\rfloor+1}-X_{\lfloor
nt\rfloor}|>\tilde{\eps} \sqrt{n}].
\end{align}
By~(\ref{UPBDRT2}), (\ref{URT12}), Lemma~\ref{Aen} and 
Condition~K, the second term of~(\ref{Finiteconvdist}) tends to~$0$
as $n\to \infty$. Hence, assuming that the following limits exist, 
we deduce that
\begin{equation}
\label{Finite altern}
\lim_{n\to \infty} \Po\Big[Z^n_{\frac{\lfloor nt\rfloor}{n}}\leq
x-\tilde{\eps} \mid \Lambda_n\Big]\leq \lim_{n\to \infty}\Po[Z^n_t\leq x \mid
\Lambda_n]\leq 
\lim_{n\to \infty} \Po\Big[Z^n_{\frac{\lfloor nt\rfloor}{n}}\leq
x+\tilde{\eps} \mid \Lambda_n\Big]
\end{equation}
for all $\tilde{\eps}>0$.
Now, suppose that we have for all $x\geq 0$ and $0<t<1$,
\begin{equation}
\label{Finite replace}
\lim_{n \to \infty} \Po\Big[Z^n_{\frac{\lfloor nt\rfloor}{n}}\leq
x\mid \Lambda_n\Big]=\int_0^x q(t,y)dy.
\end{equation}
Combining~(\ref{Finite altern}) and~(\ref{Finite replace}), we
obtain~(\ref{FFD1}) since the limit distribution $q(t,x)$ is
absolutely continuous. Our goal is now to 
show~(\ref{Finite replace}). For this, observe that
\begin{align}
\label{findimdist4}
\Po\Big[Z^n_{\frac{\lfloor nt\rfloor}{n}}\leq x\mid \Lambda_n\Big]
 &=\frac{1}{\Po[\Lambda_n]}\int_0^{\frac{xn^{1/2}}{\lfloor
nt\rfloor^{1/2}}}\Po[Z^{\lfloor nt\rfloor}_1\in dy, \Lambda_{\lfloor
nt\rfloor},X_k>0,\lfloor nt\rfloor<k\leq n]\nonumber\\
 &=\frac{\Po[\Lambda_{\lfloor
nt\rfloor}]}{\Po[\Lambda_n]}\int_0^{\frac{xn^{1/2}}{\lfloor
nt\rfloor^{1/2}}}\Po^{y\sqrt{\lfloor nt \rfloor}}
\Big[Z^n_s>0, 0\leq s \leq
1-\frac{\lfloor nt\rfloor}{n}\Big]
\Po[Z^{\lfloor nt\rfloor}_1\in
dy\mid \Lambda_{\lfloor nt\rfloor}].
\end{align}
By~(\ref{UPBDRT2}), (\ref{RUIT1}), (\ref{URT12}), and~(\ref{URT26})
we have
\begin{equation}
\label{Premlemrt}
 \lim_{n\to \infty}
\frac{\Po[\Lambda_{\lfloor nt\rfloor}]}{\Po[\Lambda_n]}=t^{-\frac{1}{2}}.
\end{equation}
Using part (v) of the UCLT and Dini's theorem on uniform 
convergence of non-decreasing sequences of continuous functions, we obtain
\[
\lim_{n\to \infty}
\Po^{z\sqrt{\lfloor nt \rfloor}}\Big[Z^n_s>0, 0\leq s \leq
1-\frac{\lfloor nt\rfloor}{n}\Big]=P^z\Big[\min_{s\in [0, 1-t]}W_s>0\Big]=P[|W_{1-t}|< z]={\tilde
N}\Big(z\Big(\frac{t}{1-t}\Big)^{\frac{1}{2}}\Big)
\]
uniformly in~$z$ on every compact set of $\R_+$.
By Proposition~\ref{propcltord}, we have 
\[
\lim_{n \to \infty}
\Po[Z^{\lfloor nt\rfloor}_1\leq x\mid \Lambda_{\lfloor
nt\rfloor}]=\int_0^x y\,e^{-\frac{y^2}{2}}dy.
\]
Now, applying Lemma 2.18 of~\cite{Igle} to (\ref{findimdist4}), 
we obtain 
\begin{equation*}
\lim_{n \to \infty}\Po[Z^n_t\leq x\mid \Lambda_n]=\lim_{n \to \infty}\Po\Big[Z^n_{\frac{\lfloor nt\rfloor}{n}}\leq x\mid \Lambda_n\Big]
= \int_0^{xt^{-\frac{1}{2}}}t^{-\frac{1}{2}}
{\tilde N}\Big(y\Big(\frac{t}{1-t}\Big)^{\frac{1}{2}}\Big)ye^{-\frac{y^2}{2}} dy.
\end{equation*}
 Finally, make the change of variables $u=t^{\frac{1}{2}}y$ to obtain the
desired result.
\qed

We can now use Propositions \ref{propcltord} and \ref{prop_Z<x} to easily check that (\ref{ORTO1}), (\ref{ORTO2}) and (\ref{ORTO3}) of Theorems \ref{Theodu1} and \ref{Theodu2} are satisfied.  
This ends the proof of Theorem \ref{Theocond}.
\qed


\section{Proof of Corollary~\ref{Corro}}
\label{s_corollary}
In this last part, for the sake of brevity,
 we will use the same notation for a real number~$x$ and its
integer part~$\lfloor x \rfloor$. The interpretation of the notation
should be clear by the context where it is used. We also suppose
without loss of generality that $\sig=1$. Let us first introduce some
spaces needed in the proof of Corollary~\ref{Corro}.

 For any $l>0$, let $C_0([0,l])$ the space of continuous
functions~$f$
from $[0,l]$ into~$\R$ such that $f(0)=0$. We endow this space with
the metric 
\[
d(f,g)=\sup_{x\in [0,l]}|f(x)-g(x)|
\]
and the Borel sigma-field on $C_0([0,l])$
corresponding to the metric~$d$.
 
Then, let $C_0(\R_+)$ the space of continuous
functions~$f:\R_+\to\R$
such that $f(0)=0$. We endow this space with the metric 
\[
\mathtt{d}(f,g)=\sum_{n=1}^{\infty}2^{-n+1}\min\{1, \sup_{x\in
[0,n]}|f(x)-g(x)|\}
\]
and the Borel sigma-field on $C_0(\R_+)$
corresponding to the metric $\mathtt{d}$.  
 Next, let~$G$ be the set of functions of $C_0(\R_+)$ for which
there exists $x_0$ (depending on~$f$) such that $f(x_0)=1$.
Let us also define the set~$H$ as the set of functions of
$C_0(\R_+)$ such that there exists
$x_1=x_1(f)=\min\{s>0:f(s)=1\}$ and $f(x)=1$ for all $x\geq
x_1$; observe that~$G$ and~$H$ are closed subsets of~$C_0(\R_+)$.
 We define the continuous map~$\Psi:G\to H$
by 
\[
\Psi(f)(x)= \left\{
    \begin{array}{ll}
        f(x) & \mbox{for}~x\leq x_1,\\
        1 & \mbox{for}~x> x_1.\\
    \end{array}
\right.
\]
Now, Corollary~\ref{Corro} can be restated as follows:
%
under the conditions of Theorem~\ref{Theocond}, we have $\IP$-a.s.,
for all measurable $A\subset H$ such that $P[B_3(\cdot
\wedge
\varrho_1)\in \partial A]=0$ and all $a\geq 0$,
\begin{equation}
\label{Corroeq}
\lim_{n\to \infty}\Po[Y^n \in A, T_n \leq a \mid \Lambda'_n]
= P[B_3(\cdot \wedge \varrho_1)\in A,  \varrho_1 \leq
a].
\end{equation}
Before proving this last statement, 
let us start by denoting $R=\{Y^n\in A\}$.
 We will bound the term $\Po[R, T_n \leq a \mid \Lambda'_n]$ from
above and below, for sufficiently large~$n$.

We start with the upper bound. Let $M>0$ be an integer and $I_M=[n,
n+M]$. We obtain
\begin{align}
\label{Triden}
\Po[R, T_n\leq a \mid \Lambda'_n]
 &=\frac{1}{\Po[\Lambda'_n]}\Po[R, T_n \leq a,\Lambda'_n]\nonumber\\
 &=\frac{1}{\Po[\Lambda'_n]}\Big(\Po[R, T_n\leq a, \Lambda'_n,
X_{\hat{\tau}_n}\in I_M]
+\Po[R, T_n\leq a, \Lambda'_n,
X_{\hat{\tau}_n}\notin I_M]\Big)\nonumber\\
 &\leq \frac{1}{\Po[\Lambda'_n]}\Po[R, T_n\leq a, \Lambda'_n,
X_{\hat{\tau}_n}\in I_M]+\Po[X_{\hat{\tau}_n}\notin I_M\mid 
\Lambda'_n]
\end{align}
for all sufficiently large $n$.
The second term of the right-hand side of~(\ref{Triden}) can be
treated easily. Indeed, by the same method we used to prove
Lemma~\ref{SNLcond}, we can show that, $\IP$-a.s., for all $\eta>0$,
there exists $M>0$ such that 
\begin{equation}
\label{IRKO1}
\Po[X_{\hat{\tau}_n}\notin I_M\mid  \Lambda'_n]\leq \eta
\end{equation}
for all $n\geq 1$.
 Let $c>2a$ and observe that $R\cap \{T_n\leq a\} \in
\mathcal{F}_{\hat{\tau}_n}$, where $\mathcal{F}_{\hat{\tau}_n}$ is
the sigma-field generated by $X$ until time $\hat{\tau}_n$. For the
first term of the right-hand side of~(\ref{Triden}), we have by the
Markov property 
\begin{align*}
\frac{1}{\Po[\Lambda'_n]}\Po[R, T_n\leq a, \Lambda'_n,
X_{\hat{\tau}_n}\in I_M]
 &=\sum_{u=0}^{M}\frac{1}{\Po[\Lambda'_n]} \Po[R, T_n\leq a,
\Lambda'_n, X_{\hat{\tau}_n}= n+u]\nonumber\\
 &=\sum_{u=0}^{M}\frac{\Po^{n+u}[\Lambda_{(c-a)n^2}]}{\Po[\Lambda'_n]
\Po^{n+u}[\Lambda_{(c-a)n^2}]} \Po[R, T_n\leq
a,\Lambda'_n,X_{\hat{\tau}_n}= n+u]\nonumber\\
 &\leq
\sum_{u=0}^{M}\frac{1}{\Po[\Lambda'_n]\Po^{n+u}[\Lambda_{(c-a)n^2}]}
\Po[R, T_n\leq a,\Lambda_{(c-a)n^2} ,X_{\hat{\tau}_n}= n+u].
\end{align*}
Next,
let us define the event $E=\{X_1>0,\dots,
X_{\hat{\tau}_n}>0,\dots,X_{\hat{\tau}_n+(c-a)n^2}>0\}$.
Using the Markov property, we can write
\begin{align*}
\Po[E]
 &\leq \sum_{v=0}^{M}\Po[X_1>0,\dots,
X_{\hat{\tau}_n}=n+v,\dots,X_{\hat{\tau}_n+(c-a)n^2}>0]
+\Po[X_{\hat{\tau}_n}\notin I_M,\Lambda'_n]\nonumber\\
 &= \sum_{v=0}^{M}\Po[\Lambda'_n, 
X_{\hat{\tau}_n}=n+v]\Po^{n+v}[\Lambda_{(c-a)n^2}]+\Po[X_{\hat{\tau}
_n}\notin I_M\mid \Lambda'_n]\Po[\Lambda'_n].
\end{align*}
But, by the UCLT,
we have for all $\eps>0$ that uniformly in $v\in [0,M]$ 
and $u\in [0,M]$,
\begin{equation}
\label{ONC}
 \Big|\Po^{n+v}[\Lambda_{(c-a)n^2}]-\Po^{n+u}[\Lambda_{(c-a)n^2}]
\Big|\leq \eps
\end{equation}
for all~$n$ sufficiently large.
Therefore, we obtain for all $u\in [0,M]$,
\begin{equation}
\label{Temp11}
\Po^{n+u}[\Lambda_{(c-a)n^2}]\Po[\Lambda'_n]\geq
\Po[E]-(\eps+\eta)\Po[\Lambda'_n]
\end{equation}
for all~$n$ sufficiently large.
Now, let us bound the first term of the right-hand side
of~(\ref{Temp11}) from below.
Fix some $\delta>0$. We write
\begin{align}
\label{UU2}
\Po[E]
&\geq \Po[E,\hat{\tau}_n\leq(a+\delta)n^2]\nonumber\\
 &\geq\Po[\Lambda_{((c+\delta)n^2+3)},\hat{\tau}_n\leq(a+\delta)n^2]
\nonumber\\
 &\geq
\Po[\Lambda_{((c+\delta)n^2+3)}]\Po[\hat{\tau}_n\leq(a+\delta)n^2\mid
\Lambda_{((c+\delta)n^2+3)}].
\end{align}
Finally, by~(\ref{Triden}), (\ref{IRKO1}), (\ref{Temp11}) 
and~(\ref{UU2}) we obtain $\IP$-a.s.,
\begin{align}
\label{RUTU}
\Po[R, T_n\leq a \mid \Lambda'_n]
 &\leq 
\frac{(\Po[\Lambda_{((c+\delta)n^2+3)}])^{-1}\Po[\Lambda_{(c-a)n^2}]
\Po[R, T_n\leq
a\mid\Lambda_{(c-a)n^2}]}{\Po[\hat{\tau}_n\leq(a+\delta)n^2\mid
\Lambda_{((c+\delta)n^2+3)}]-(\eps+\eta)
\Po[\Lambda'_n](\Po[\Lambda_{((c+\delta)n^2+3)}))^{-1}}+\eta 
\end{align}
for all sufficiently large~$n$.

We now estimate the term $\Po[R, T_n\leq a \mid \Lambda'_n]$ from
below.
Let us write
\begin{equation}
 \label{Triden2}
\Po[R, T_n\leq a \mid \Lambda'_n]
=\frac{1}{\Po[\Lambda'_n]}\Po[R, T_n\leq a,\Lambda'_n]
\geq \frac{1}{\Po[\Lambda'_n]}\Po[R, T_n\leq a, \Lambda'_n,
X_{\hat{\tau}_n}\in I_M].
\end{equation}
Then, we have by the Markov property 
\begin{align}
\frac{1}{\Po[\Lambda'_n]}\Po[R, T_n\leq a, \Lambda'_n,
X_{\hat{\tau}_n}\in I_M]
 &=\sum_{u=0}^{M}\frac{1}{\Po[\Lambda'_n]} \Po[R, T_n\leq a,
\Lambda'_n, X_{\hat{\tau}_n}= n+u]\nonumber\\
 &=\sum_{u=0}^{M}\frac{\Po^{n+u}[\Lambda_{(c-a)n^2}]}{\Po[\Lambda'_n]
\Po^{n+u}[\Lambda_{(c-a)n^2}]} \Po[R, T_n\leq
a,\Lambda'_n,X_{\hat{\tau}_n}= n+u]\nonumber\\
 &\geq
\sum_{u=0}^{M}\frac{1}{\Po[\Lambda'_n]\Po^{n+u}[\Lambda_{(c-a)n^2}]}
\Po[R, T_n\leq a,\Lambda_{cn^2} ,X_{\hat{\tau}_n}= n+u].
\label{Foufur}
\end{align}
Again using the Markov property, we can write
\begin{align*}
\Po[E]
 &\geq \sum_{v=0}^{M}\Po[X_1>0,\dots,
X_{\hat{\tau}_n}=n+v,\dots,X_{\hat{\tau}_n+(c-a)n^2}>0]\nonumber\\
 &\geq \sum_{v=0}^{M}\Po[\Lambda'_n, 
X_{\hat{\tau}_n}=n+v]\Po^{n+v}[\Lambda_{(c-a)n^2}].
\end{align*}
Using~(\ref{ONC}), we obtain for all $u\in [0,M]$,
\begin{equation}
\label{Temp12}
\Po^{n+u}[\Lambda_{(c-a)n^2}]\Po[\Lambda'_n]\leq
\Po[E]+\eps\Po[\Lambda'_n]
\end{equation}
for all sufficiently large~$n$.
Then, as $\hat{\tau}_n\geq 1$, we have
\begin{equation}
\label{UU3}
\Po[E]\leq \Po[\Lambda_{(c-a)n^2}].
\end{equation}
Finally, by~(\ref{Triden2}), (\ref{Foufur}), (\ref{Temp12}) 
and~(\ref{UU3}), we
obtain $\IP$-a.s.,
\begin{align}
\label{RUTU2}
\Po[R, T_n\leq a \mid \Lambda'_n]
&\geq  \frac{(\Po[\Lambda_{(c-a)n^2}])^{-1}\Po[\Lambda_{cn^2}]\Po[R,
T_n\leq a\mid\Lambda_{cn^2}]}{1+\eps
\Po[\Lambda'_n](\Po[\Lambda_{(c-a)n^2}])^{-1}}
\end{align}
for all sufficiently large~$n$.

Our intention is now to take the $\limsup$ as $n\to \infty$ 
in~(\ref{RUTU}). Before this, observe that by~(\ref{UPBDRT2}), 
(\ref{RUIT1}), (\ref{URT12}) and~(\ref{URT26}) we have for $\eps\leq 1$,
\begin{align}
\lim_{n\to \infty}
\frac{\Po[\Lambda_{(c-a)n^2}]}{\Po[\Lambda_{((c+\delta)n^2+3)}]}
&=\sqrt{\frac{c+\delta}{c-a}}, \label{Xant1}\\
\limsup_{n\to \infty} 
\frac{\Po[\Lambda'_n]}{\Po[\Lambda_{((c+\delta)n^2+3)}]}&\leq \limsup_{n\to \infty} 
\frac{\Po[A_{\eps,n^2}]}{\Po[\Lambda_{((c+\delta)n^2+3)}]}
\leq \gamma_1\sqrt{c+\delta} \label{Xant2}
\end{align}
for some constant $\gamma_1$.
By the usual scaling, from the Brownian meander $W^+$ on $[0,1]$ it is possible to define the 
Brownian meander~$W^+_t$ on any finite interval $[0,t]$: $W^+_t(\cdot):=\sqrt{t}W^+(\cdot/t)$. 
Thus, Theorem~\ref{Theocond} implies that
\begin{equation}
\label{Xant3}
\lim_{n\to \infty}
\Po[\hat{\tau}_n\leq(a+\delta)n^2\mid
\Lambda_{((c+\delta)n^2+3)}]=P\Big[\sup_{0\leq s\leq
(a+\delta)}W^+_{c+\delta}(s)\geq 1\Big].
\end{equation}

Denoting 
by $\mathcal{U}_a$ the measurable set of functions~$f$ in~$H$ such
that $f(a)=1$ and by $\pi_l$ the projection map from $C_0(\R_+)$
onto $C_0([0,l])$, we have
\begin{align*}
\Po[R, T_n\leq a\mid\Lambda_{(c-a)n^2}]
 &=\Po[Z^{n^2}_{\cdot \wedge T_n}\in A \cap
\mathcal{U}_a\mid\Lambda_{(c-a)n^2}]\nonumber\\
 &=\Po[Z^{n^2}\in \Psi^{-1}(A \cap
\mathcal{U}_a)\mid\Lambda_{(c-a)n^2}]\nonumber\\
 &=\Po[Z^{n^2}_{\cdot \wedge (c-a)}\in \pi_{c-a}(\Psi^{-1}(A \cap
\mathcal{U}_a))\mid\Lambda_{(c-a)n^2}].
\end{align*}
The next step is to show that 
\begin{align}
\label{Xant4}
\lim_{n\to \infty}\Po[Z^{n^2}_{\cdot \wedge (c-a)}\in
\pi_{c-a}(\Psi^{-1}(A \cap
\mathcal{U}_a))\mid\Lambda_{(c-a)n^2}]
&=P[W^+_{c-a}\in \pi_{c-a}(\Psi^{-1}(A \cap \mathcal{U}_a))]
\end{align}
where $W^+_{c-a}$ is the Brownian meander on $[0,c-a]$.
 As the law of the Brownian meander on $[0,c-a]$ is absolutely
continuous with respect to the law of the three dimensional Bessel
process $B_3$ on $[0,c-a]$ (see~\cite{Imhof} section~4), to
prove~(\ref{Xant4}) we will show that 
\begin{equation}
 P[B_3(\cdot \wedge c-a)\in \partial \{\pi_{c-a}(\Psi^{-1}(A \cap
\mathcal{U}_a))\}]=0.
\end{equation}
Observe that, as $\pi_{c-a}$ is a projection, we have 
\begin{align*}
P[B_3(\cdot \wedge c-a)\in \partial \{\pi_{c-a}(\Psi^{-1}(A
\cap \mathcal{U}_a))\}]
 &\leq P[B_3(\cdot \wedge c-a)\in \pi_{c-a} \partial \{\Psi^{-1}(A
\cap \mathcal{U}_a)\}]\nonumber\\
 &= P[B_3\in \partial \{\Psi^{-1}(A \cap \mathcal{U}_a)\}].
\end{align*}
Now, as $\Psi$ is a continuous map, we have
\begin{align}
P[B_3\in \partial \{\Psi^{-1}(A \cap \mathcal{U}_a)\}]
 &\leq P[B_3\in \Psi^{-1}(\partial\{A \cap
\mathcal{U}_a\})]\nonumber\\
 &\leq P[B_3\in \Psi^{-1}(\partial A \cup \partial
\mathcal{U}_a)]\nonumber\\
 &\leq P[B_3(\cdot \wedge \varrho_1)\in \partial A]+P[\varrho_1=a].
\end{align}
 By hypothesis, $P[B_3(\cdot \wedge \varrho_1)\in \partial A]=0$. As
the law of~$\varrho_1$ is absolutely continuous with respect to the
Lebesgue measure (see~\cite{Imhof} Theorem~4), we also have
$P[\varrho_1=a]=0$. This proves~(\ref{Xant4}).

Then, we want to take the $\liminf$ as $n\to \infty$
in~(\ref{RUTU2}). Before this, notice that
\begin{align}
 \lim_{n\to
\infty}\frac{\Po[\Lambda_{cn^2}]}{\Po[\Lambda_{(c-a)n^2}]}
&=\sqrt{\frac{c-a}{c}},\label{Xant5}\\
 \limsup_{n\to \infty}
\frac{\Po[\Lambda'_n]}{\Po[\Lambda_{(c-a)n^2}]}&\leq
\gamma_2\sqrt{c-a}\label{Xant6}
\end{align}
for some constant $\gamma_2$.
By the same argument we used to prove~(\ref{Xant4}), we have
\begin{align}
\label{Xant7}
\lim_{n\to \infty}\Po[R, T_n\leq a\mid\Lambda_{cn^2}]
= P[W^+_{c}\in \pi_{c}(\Psi^{-1}(A \cap \mathcal{U}_a))]
\end{align}
where $W^+_{c}$ is the Brownian meander on $[0,c]$.
Then, define 
$V_{l}=\{W^+_{l}\in \pi_{l}(\Psi^{-1}(A \cap \mathcal{U}_a))\}$
for $l\in \{c-a,c\}$.
Combining~(\ref{Xant1}), (\ref{Xant2}), (\ref{Xant3}),
(\ref{Xant4}), (\ref{Xant5}), (\ref{Xant6}) and~(\ref{Xant7}) we see
that 
\begin{align}
\frac{P[V_{c}]\sqrt{\frac{c-a}{c}} }{1+ \gamma_2\eps \sqrt{c-a} }
&\leq\liminf_{n\to \infty} \Po[R, T_n\leq a \mid
\Lambda'_n] \leq\limsup_{n\to \infty} \Po[R, T_n\leq a \mid
\Lambda'_n]\nonumber\\
&\leq\frac{P[V_{c-a}]\sqrt{\frac{c+\delta}{c-a}}}{P[\sup_{0\leq s\leq
(a+\delta)}W^+_{c+\delta}(s)\geq
1]-\gamma_1(\eps+\eta)\sqrt{c+\delta} }+\eta.
\end{align}

 Now, take $\eps=\eta=c^{-1}$ and $\delta=\sqrt{c}$ and let $c$ tend
to infinity. Since 
\[
 P[W_l^+\in \pi_{l}(\Psi^{-1}(A \cap \mathcal{U}_a))]=P[W_l^+(\cdot
\wedge a)\in \pi_{a}(\Psi^{-1}(A \cap \mathcal{U}_a))],
\]
we have by Lemma 11-1 of~\cite{BY}
\[
 \lim_{c\to \infty} P[V_{l}]=P[B_3(\,\cdot\wedge  \varrho_1)\in A,
\varrho_1\leq a]
\]
for $l\in \{c-a,c\}$.

 The last thing we have to check to obtain~(\ref{Corroeq}) is that 
\begin{equation}
\label{GRANFIN2}
 \lim_{c\to \infty}P\Big[\sup_{0\leq s\leq
(a+\delta)}W^+_{c+\delta}(s)< 1\Big]=0.
\end{equation}
First, we start by noting that by scaling property
\begin{align*}
 P\Big[\sup_{0\leq s\leq (a+\delta)}W^+_{c+\delta}(s)< 1\Big]
 &=P\Big[(c+\delta)^{\frac{1}{2}}\sup_{0\leq s\leq
(a+\delta)}W^+\Big(\frac{s}{c+\delta}\Big)< 1\Big]\nonumber\\
 &=P\Big[\sup_{0\leq s\leq \frac{a+\delta}{c+\delta}}W^+(s)<
(c+\delta)^{-\frac{1}{2}}\Big]\nonumber\\
 &\leq P\Big[W^+\Big(\frac{a+\delta}{c+\delta}\Big)\leq
(c+\delta)^{-\frac{1}{2}}\Big]
\end{align*}
where $W^+$ is a Brownian meander on $[0,1]$.
 This last term is easily computable using the transition density
function from $(0,0)$ of $W^+$ given in~(\ref{GRANFIN}). Let
$u=\frac{a+\delta}{c+\delta}$,
\begin{align*}
P\Big[W^+\Big(\frac{a+\delta}{c+\delta}\Big)\leq
(c+\delta)^{-\frac{1}{2}}\Big]
 &=\int_{0}^{(c+\delta)^{-\frac{1}{2}}}u^{-\frac{3}{2}}x\exp{\Big(-\frac{x^2}{2u}
\Big)}{\tilde N}(x(1-u)^{-\frac{1}{2}})dx
\end{align*}
 Let us make the change of variable $y=(c+\delta)^{\frac{1}{2}}x$ in the
right-hand side integral. Then, we obtain 
\begin{align*}
P\Big[W^+\Big(\frac{a+\delta}{c+\delta}\Big)\leq
(c+\delta)^{-\frac{1}{2}}\Big]
 &=\int_{0}^{1}\frac{(c+\delta)^{\frac{1}{2}}}{(a+\delta)^{\frac{3}{2}}}y\exp{
\Big(-\frac{y^2}{2(a+\delta)}\Big)}{\tilde N}(y(c-a)^{-\frac{1}{2}})dy.
\end{align*}
 Now, making the change of variable $z=(c-a)^{\frac{1}{2}}u$ in the following
integral
\begin{equation*}
 {\tilde
N}(y(c-a)^{-\frac{1}{2}})=\sqrt{\frac{2}{\pi}}\int_0^{y(c-a)^{-\frac{1}{2}}}\exp{\Big(-\frac{u^2}{2}\Big)}
du
\end{equation*}
we obtain
\begin{align}
\label{JUNI}
\lefteqn{
P\Big[W^+\Big(\frac{a+\delta}{c+\delta}\Big)\leq
(c+\delta)^{-\frac{1}{2}}\Big]
}\phantom{**********}\nonumber\\
 &=\Big(\frac{c+\delta}{c-a}\Big)^{\frac{1}{2}}\frac{1}{(a+\delta)^{\frac{3}{2}}}
\int_{0}^{1}\int_0^y
 y\exp{\Big(-\frac{y^2}{2(a+\delta)}\Big)}\exp{\Big(-\frac{z^2}{
2(c-a) }\Big)}dz\,dy\nonumber\\
 &\leq
\Big(\frac{c+\delta}{c-a}\Big)^{\frac{1}{2}}\frac{1}{(a+\delta)^{\frac{3}{2}}}\int_{0
}^{1}\int_0^y y\,dz\,dy\nonumber\\
 &=\frac{1}{3}\Big(\frac{c+\delta}{c-a}\Big)^{\frac{1}{2}}\frac{1}{
(a+\delta)^{\frac{3}{2}}}.
\end{align}
 Taking $\delta=\sqrt{c}$ and letting $c\to \infty$ in~(\ref{JUNI}),
we obtain~(\ref{GRANFIN2}). This concludes the proof
of~(\ref{Corroeq}).
\qed

\section*{Acknowledgments}
C.G.\ is grateful to FAPESP (grant 2009/51139--3) for financial
support. S.P.\ was partially supported by
CNPq (grant 300886/2008--0). 
Both also thank CNPq (472431/2009--9) and FAPESP (2009/52379--8)  
for financial support. Both authors thank the anonymous referee for valuable comments which 
allowed to improve the first version of this paper.

\end{document}